\documentclass[preprint,number]{elsarticle}
\biboptions{sort&compress}
\pdfoutput=1

\usepackage{graphicx}
\usepackage{longtable}
\usepackage{multirow}
\usepackage{algorithm}

\usepackage{subfigure}

\usepackage[draft=false,pagebackref=false]{hyperref}

\hypersetup{
    linktoc=page,
    colorlinks=true,        
    linkcolor=black,          
    citecolor=black,          
    filecolor=black,         
    urlcolor=black           
}

\usepackage{amsmath,amssymb}

\graphicspath{{./}{./figs/}}

  \def\clap#1{\hbox to 0pt{\hss#1\hss}}

\usepackage{bm}
\providecommand{\mat}[1]{\bm{#1}}%
\renewcommand{\vec}[1]{\mathbf{#1}}
\newcommand{\vecalt}[1]{\bm{#1}}

\providecommand{\mA}{\ensuremath{\mat{A}}}

\providecommand{\mC}{\ensuremath{\mat{C}}}

\providecommand{\mW}{\ensuremath{\mat{W}}}

\providecommand{\vf}{\ensuremath{\vec{f}}}

\providecommand{\vu}{\ensuremath{\vec{u}}}

\providecommand{\vw}{\ensuremath{\vec{w}}}
\providecommand{\vx}{\ensuremath{\vec{x}}}

\providecommand{\vpi}{\ensuremath{\vecalt{\pi}}}

\newcommand{\sS}{\mathcal{S}}

\newcommand{\bmat}[1]{\begin{bmatrix}#1\end{bmatrix}}

\newcommand{\argmax}[1]{\underset{#1}{\mathrm{argmax}}}
\newcommand{\argmin}[1]{\underset{#1}{\mathrm{argmin}}}

\newcommand{\st}{\text{subject to }}

\newcommand{\ymax}{y_{\text{max}}}

\usepackage{memory-macros}

\renewenvironment{algorithm}[
1][\relax]{\refstepcounter{algorithm}%
\addcontentsline{loa}{algorithm}%
    {\protect\numberline{Algorithm~\thealgorithm}{\ignorespaces#1}}%
\par\vspace{1\baselineskip}%
\expandafter\ifx#1\relax
\parindent0pt {\scshape\bfseries Algorithm~\thealgorithm.}\\
\else
\parindent0pt {\scshape\bfseries
Algorithm~\thealgorithm.}\enspace{\bfseries#1.}\\
\fi}
{\vspace{1\baselineskip}\par}

\begin{document}


\title{Exploiting Active Subspaces to Quantify Uncertainty in the Numerical Simulation of the {H}y{S}hot {II} Scramjet}

\author{P.G.~Constantine\corref{cor}\fnref{csm}}
\ead{paul.constantine@mines.edu}
\ead[url]{http://inside.mines.edu/~pconstan}

\author[ctr]{M.~Emory\fnref{ctr}}
\ead{memory@stanford.edu}

\author[umd]{J.~Larsson\fnref{umd}}
\ead{jola@umd.edu}

\author[ctr]{G.~Iaccarino\fnref{ctr}}
\ead{jops@stanford.edu}

\address[csm]{Department of Applied Mathematics and Statistics, 1500 Illinois Street, Colorado School of Mines, Golden, Colorado, 80401}
\address[ctr]{Center for Turbulence Research, 488 Escondido Mall, Building 500, Stanford University, Stanford, California, 94305}
\address[umd]{Department of Mechanical Engineering, 3149 Glenn L.~Martin Hall, University of Maryland, College Park, MD, 20742}

\cortext[cor]{Corresponding author}

\begin{abstract}
We present a computational analysis of the reactive flow in a hypersonic scramjet engine with focus on effects of uncertainties in the operating conditions. We employ a novel methodology based on \emph{active subspaces} to characterize the effects of the input uncertainty on the scramjet performance. The active subspace identifies one-dimensional structure in the map from simulation inputs to quantity of interest that allows us to reparameterize the operating conditions; instead of seven physical parameters, we can use a single derived active variable. This dimension reduction enables otherwise infeasible uncertainty quantification, considering the simulation cost of roughly 9500 CPU-hours per run. For two values of the fuel injection rate, we use a total of 68 simulations to (i) identify the parameters that contribute the most to the variation in the output quantity of interest, (ii) estimate upper and lower bounds on the quantity of interest, (iii) classify sets of operating conditions as safe or unsafe corresponding to a threshold on the output quantity of interest, and (iv) estimate a cumulative distribution function for the quantity of interest. 
\end{abstract}


\maketitle

\section{Introduction}
\label{sec:intro}


\noindent Over the past decade there has been a renewed interest in numerical simulations of high-speed air-breathing propulsion systems for hypersonic vehicles, driven by fundamental advancements in simulation tools for high Mach-number flight conditions. These tools are essential in the design process of such vehicles due to the challenges and costs associated with traditional physical prototyping. In light of the limited experience with sustained hypersonic flight, numerical simulations can be of critical value to understanding the behavior of this system, in particular characterizing the safe operability limits of the propulsion system.

\begin{figure}[htb!]
  \includegraphics[width=1\linewidth]{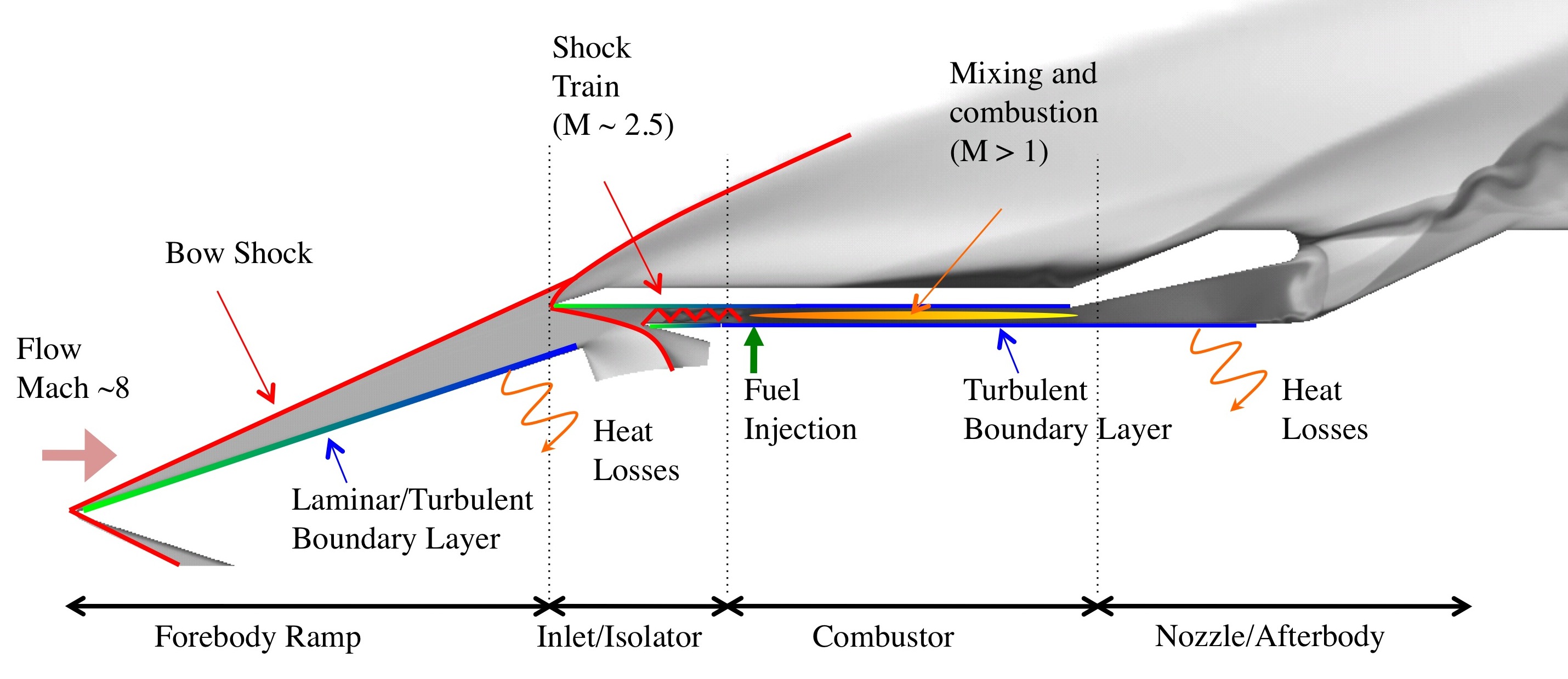}
  \caption{Side view of the HyShot II scramjet outlining the various physics within the different parts of the geometry. Contours show numerical Schlieren results from a 2D simulation.}
  \label{fig:hyshot-phys-diagram}
\end{figure}

Supersonic combustion engines (scramjets) are an economic alternative to rockets because they do not require on-board storage of the oxidizer. The HyShot II scramjet, shown in Figure \ref{fig:hyshot-phys-diagram}, was designed to demonstrate supersonic combustion during flight in a simple configuration~\cite{Hass2005,Smart2006}, and it has since been the subject of multiple ground-based experimental campaigns in the High Enthalpy shock tunnel G{\"o}ttingen (HEG) of the German Aerospace Center (DLR)~\cite{Gardner,Gardner2004,Schramm2008,Hannemann2009,Hannemann2010,laurence:11,laurence:12,laurence:13,laurence:14}.
The simplicity of the configuration and the availability of
experimental data has made the HyShot II scramjet the subject of
multiple computational investigations, based on either 
Reynolds-averaged Navier-Stokes (RANS) simulations~\cite{Pecnik2012}
or large-eddy simulations (LES)~\cite{Fureby2010,Chapius2013,larsson:14}.

Although the geometry depicted in Figure \ref{fig:hyshot-phys-diagram} (and in detail in Figure \ref{fig:hyshot-geom}) is simple and without moving parts, accurately predicting the internal flow structure is a very challenging problem. The difficulty results from the physical phenomena encountered within the engine when traveling at high velocity---namely turbulence, shocks, boundary layers, mixing, and combustion---all of which must be modeled. Required closures and disparate spatial and temporal scales of the various physical processes cause these simulations to be fairly expensive.


Scramjet operation requires a careful balance between maximizing thrust and maintaining stable and safe operation. At the conditions of interest here, the flow inside the HyShot II scramjet combustor has three different regimes \cite{laurence:14,larsson:14}.
At low fuel-air equivalence ratios (about $\lesssim 0.38$), the flow is supersonic throughout. At higher fuel-air equivalence ratios, a stable flow with a shock-train in the combustor develops.
Finally, at sufficiently high fuel-air equivalence ratios, the flow will unstart---a potentially catastrophic failure mode associated with a large loss in thrust. Scramjet design and operation is complicated by the fact that the maximum thrust occurs very close to the boundaries between these different regimes.

Previous efforts \cite{Waltrup1973,Wagner2009,Byrne2000} studied various aspects of this problem both experimentally and computationally, but they largely ignored the effect of the uncertainties on the estimation of the operability limit. These efforts enhance understanding of the physics but do not provide sufficient confidence in the quantitative estimates of the maximum fuel flow rate compatible with safe operations. As a consequence safety margins based on experience are applied to the estimate.
To produce more trustworthy uncertainty estimates, it is necessary to explicitly identify, estimate, and account for uncertainties that affect the operability limit. Many sources of uncertainty are present in the scenario of interest; both variability in the operating environment of the vehicle (e.g., fluctuations in temperature) and inadequacy of the constitutive models (e.g., the combustion model) introduce uncertainty in the simulation result.
Our prior work studies the effects of uncertainties in the RANS turbulence models and how to account for them in a physically meaningful manner \citep{Emory2013}. We have also compared the effects of uncertainties in the combustion chemistry to effects of uncertain operating conditions for HyShot II, and we found that the effects of operating conditions were much greater than the effects of the chemistry parameters~\cite{icossar2013}. The current paper focuses on the HyShot II operating environment vehicle studied by DLR in the HEG shock tunnel \cite{Schramm2008,Gardner2004,Gardner}. In this context, the main source of variability corresponds to incomplete knowledge and controllability of the conditions within the HEG shock tunnel. A more comprehensive study that includes the effects of the operating condition variability and the uncertainties introduced by the physical models is currently ongoing and will be presented in a forthcoming paper.


Monte Carlo and related random sampling methods are commonly used to quantify uncertainty in a simulation's predictions---e.g., with moments, quantiles, or empirical density functions---given uncertainty in its inputs~\cite{mcbook,smith2013uncertainty}. Estimates are unbiased, and one can develop confidence intervals for the uncertainty analysis from the Central Limit Theorem. However, when the simulation is expensive, sampling can be impractical for accurate statistics due to the slow $\mathcal{O}(n^{-1/2})$ convergence rate (where $n$ is the number of samples). An alternative approach for rapidly converging, biased estimation employs response surfaces---where one uses a few carefully selected runs of the expensive model to construct a cheaper response surface that is subsequently sampled. Popular response surface models in uncertainty quantification include polynomial chaos~\cite{ghanem1991stochastic,xiu2002wiener}, stochastic collocation~\cite{xiu2005high}, and Gaussian process regression~\cite{rasmussen2006gaussian,Koehler1996}. For polynomial variants, moments reduce to numerical quadrature on the integral formulation of averages. Response surfaces are most appropriate when the number of input parameters is sufficiently small and the outputs are sufficiently smooth functions of the inputs. 

We supplement the quantity of interest's statistical measures with its range; determining the range is posed as optimization. We employ \emph{active subspaces}~\cite{constantine2013active} to discover that the quantity of interest can be well represented by a univariate function of the \emph{active variable} derived from the model's inputs. We then exploit this low-dimensional approximation to (i) estimate the cumulative distribution function of the quantity of interest, (ii) determine the range of the quantity of interest, and (iii) identify the sets of parameters corresponding to safe operation of the scramjet. As a byproduct, the active subspace reveals a subset of the model's parameters whose perturbations produce the greatest change in the quantity of interest, i.e., a sensitivity analysis. For this particular scramjet model, discovering the low-dimensional structure and quantifying uncertainty requires enough simulation runs to fit and validate a global, least-squares-fit linear approximation of the quantity of interest as a function of scramjet inputs. In practice, the cost scales linearly with the dimension of the input space---much lower than the scaling for more complicated approximation procedures. This cost is remarkably small considering both the simulation's complexity and its dependence on seven independent input parameters. 

The remainder of this paper is broadly structured in four sections: describing the HyShot II simulation in Section \ref{sec:rans}, characterizing the uncertainties in the operating conditions in Section \ref{sec:sources}, outlining the technique for discovering the active subspace in Section \ref{sec:uq}, and applying the active subspace to quantify uncertainties in Section \ref{sec:ashyshot}. We conclude with a brief summary.

\section{Methodology for HyShot II simulations}
\label{sec:rans}

\subsection{Geometry and computational grids}

\begin{figure}[h!tb]
  \includegraphics[width=\linewidth]{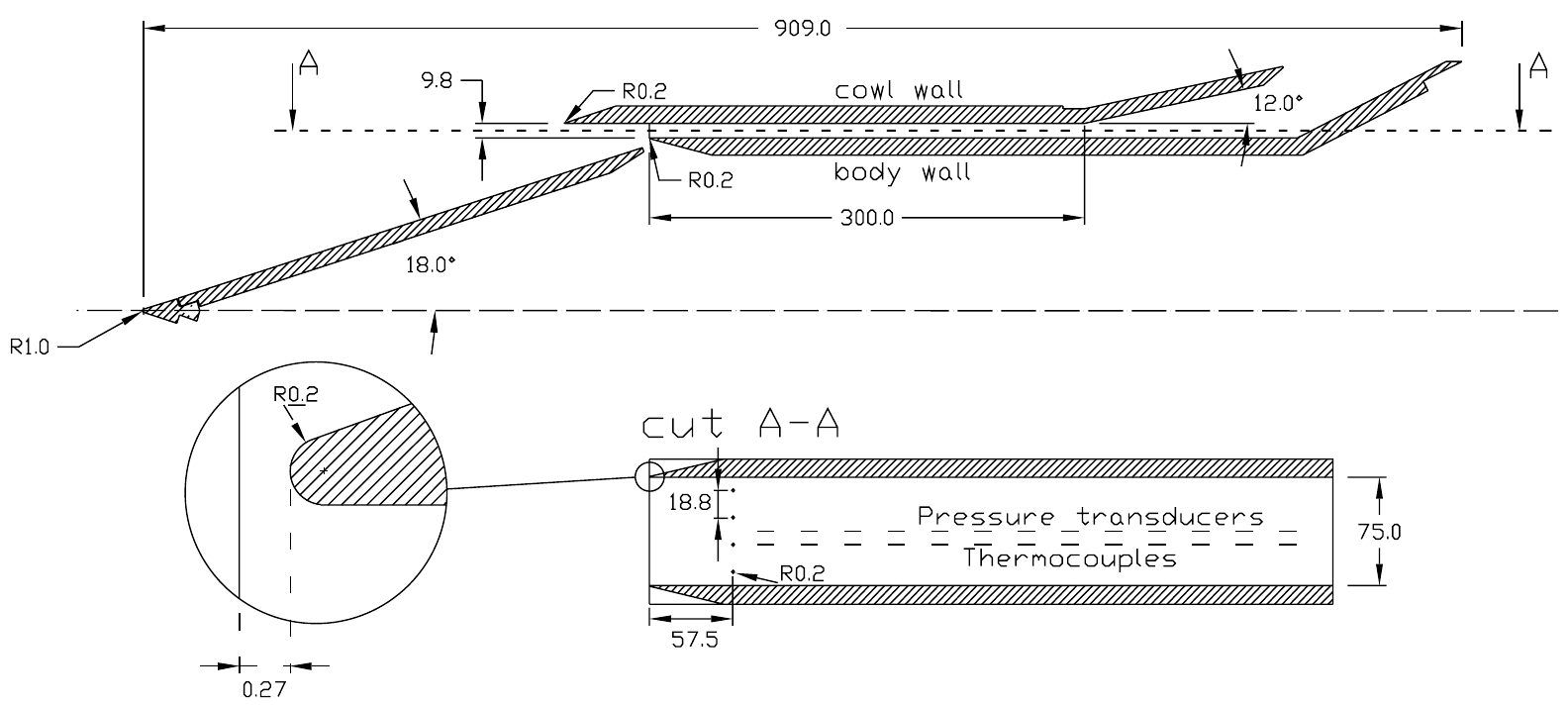}
  \caption{HyShot II scramjet geometry, taken from \cite{Pecnik2012}.}
  \label{fig:hyshot-geom}
\end{figure}

\noindent Figure \ref{fig:hyshot-geom} shows the HyShot II geometry including the intake ramp (forebody), isolator, combustion chamber, and nozzle details. The actual HEG system is much larger and includes the mounting structure and the instrumentation of the $1:1$ scale model. The geometry in Figure \ref{fig:hyshot-geom} corresponds to one-half the flight vehicle, which had symmetrically mounted fueled and un-fueled engines on top of a rocket. Both dual and single engine models have been investigated in the HEG~\cite{Schramm2008,Gardner2004,Gardner}.

The intake ramp is significantly wider than the combustion chamber, and thus the flow entering the isolator/combustor is very close to uniform in the spanwise direction. The gap between the intake ramp and body wall of the combustion chamber is a boundary-layer and isolator shock (emanating from the cowl wall leading edge) bleed channel. Within the combustor there are four fuel injection ports spanning the width of the channel, shown in the slice view \verb=A-A=. The model was mounted at a nominal angle-of-attack of $3.6^\circ$.

In order to reduce the computational cost of each simulation the HyShot II configuration is separated into two domains. The first is a 2D representation of the intake ramp and entrance to the combustion chamber, denoted as forebody and isolator in Figure \ref{fig:hyshot-phys-diagram}, respectively. The second domain is a fully 3D representation of the combustion chamber and exit nozzle, including the fuel injection ports. This domain decomposition assumes that the combustion chamber inflow conditions, i.e., the flow dynamics upstream of fuel injection, are two-dimensional (no variation in the spanwise dimension). This assertion has been verified by DLR through comparison of 2D and 3D intake ramp simulations \cite{Karl2008}. Leveraging these observations, a similar domain decomposition has been applied by a variety of researchers numerically investigating the HyShot II scramjet \cite[cf.][]{Fureby2010,Karl2008,Karl2011}.

The domain decomposition allows us to address the different modeling requirements in the 2D and 3D domains. Over the forebody and through the isolator, turbulence and transition phenomena dictate the flow structures. In contrast, the combustion chamber requires modeling of mixing and flow structures due to the injected fuel. The 2D simulations are dramatically less expensive  due to relaxed mesh resolution and modeling complexity. We next describe the meshes used in this analysis.


\subsubsection{2D forebody}
\label{sec:2d}
\begin{figure}[htb]
\begin{center}
 \includegraphics[width=0.75\linewidth]{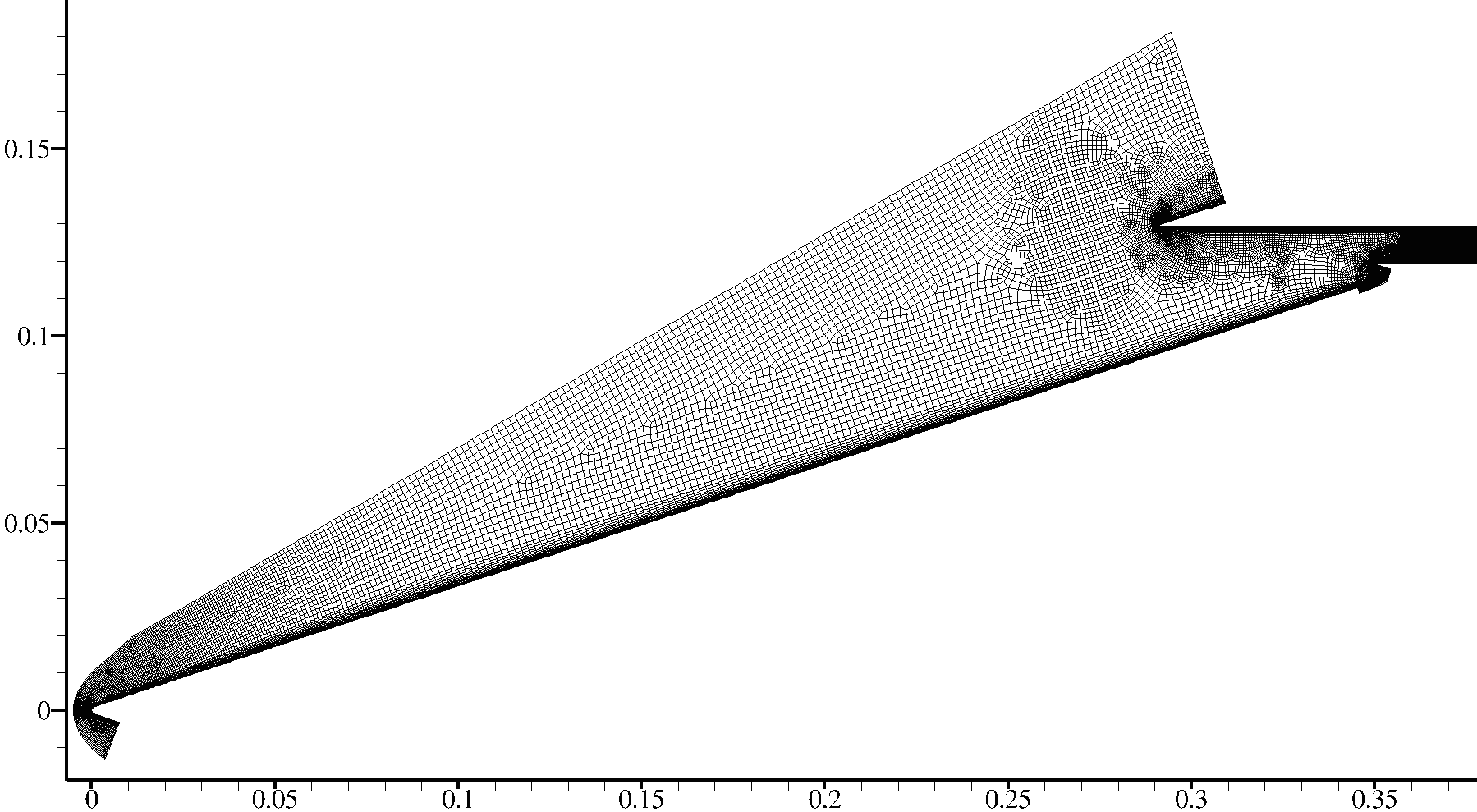}
\caption{The 2D forebody and isolator mesh used for the HyShot II simulations.}
 \label{fig:hyshot-2D-geom}
 \end{center}
\end{figure}

\noindent The grid for the 2D forebody/ramp is shown in Figure \ref{fig:hyshot-2D-geom}, where the downstream boundary is located at the fuel injection ports. In the wall-normal direction the grid is designed such that the first row of cells adjacent to the wall have $y_1^+ = 1$, resulting in just under $50$k control volumes in the domain. In regions far from solid boundaries the cells are unstructured, which reduces numerical artifacts related to poor shock-grid alignment.

\subsubsection{3D combustion chamber}
\label{sec:3d}
\noindent To reduce the cost of simulating the 3D domain we take advantage of symmetries and simulate only $1/8^{th}$ of the total combustor. In the work of Pe\v{c}nik \emph{et al.} \cite{Pecnik2010,Pecnik2012} the results using a $1/2\!-\!$span domain (including the sidewall and two injection ports) and a $1/8\!-\!$span domain (including half an injection port and using symmetry planes in the spanwise directions) were compared. Those authors concluded that while shocks emanating
from the sidewalls influence certain quantities of interest, in general the $1/8\!-\!$span domain is acceptable for use in the analysis of the scramjet. The same conclusion has been reached by DLR and other analysts, all of whom perform HyShot II simulations using a $1/8\!-\!$span domain \cite[cf.][]{Fureby2010,Karl2008,Karl2011,Chen2011}.

\begin{figure}[h!tb]
\begin{center}
 \subfigure[]{\includegraphics[width=0.44\linewidth]{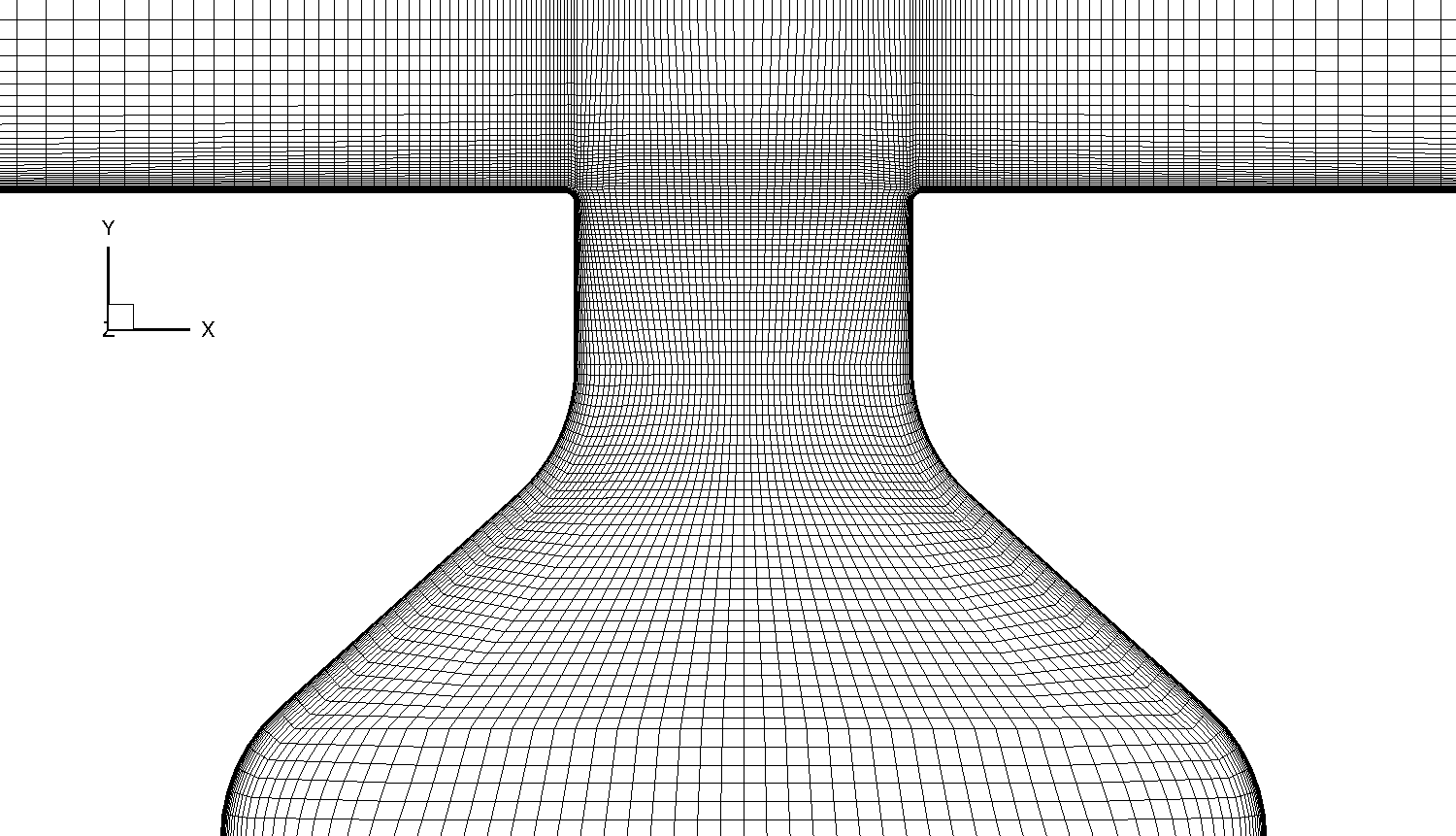}\label{fig:hyshot-3D-geom-a}}
 \subfigure[]{\includegraphics[width=0.44\linewidth]{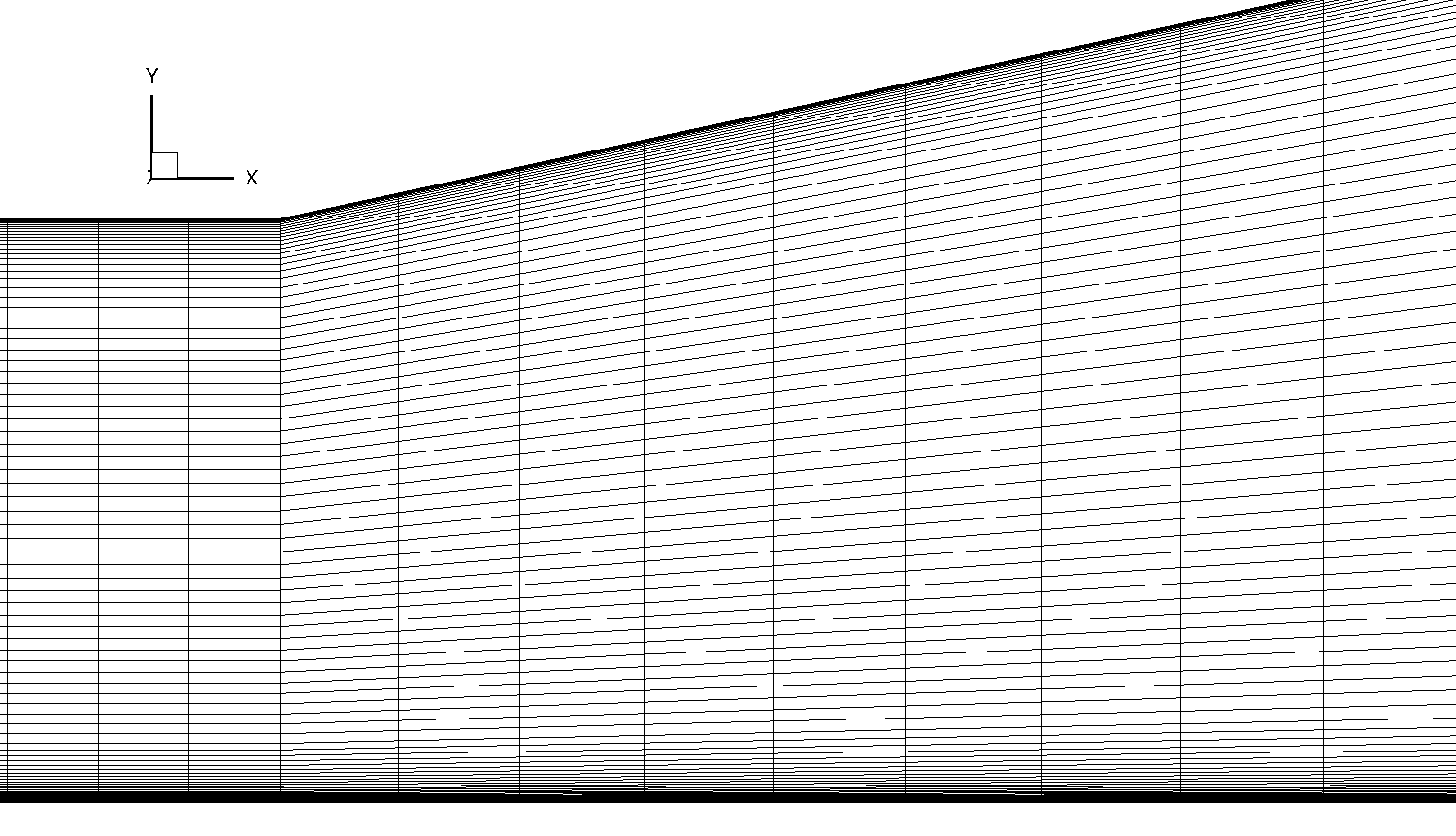}\label{fig:hyshot-3D-geom-b}}
 \subfigure[]{\includegraphics[width=.9\linewidth]{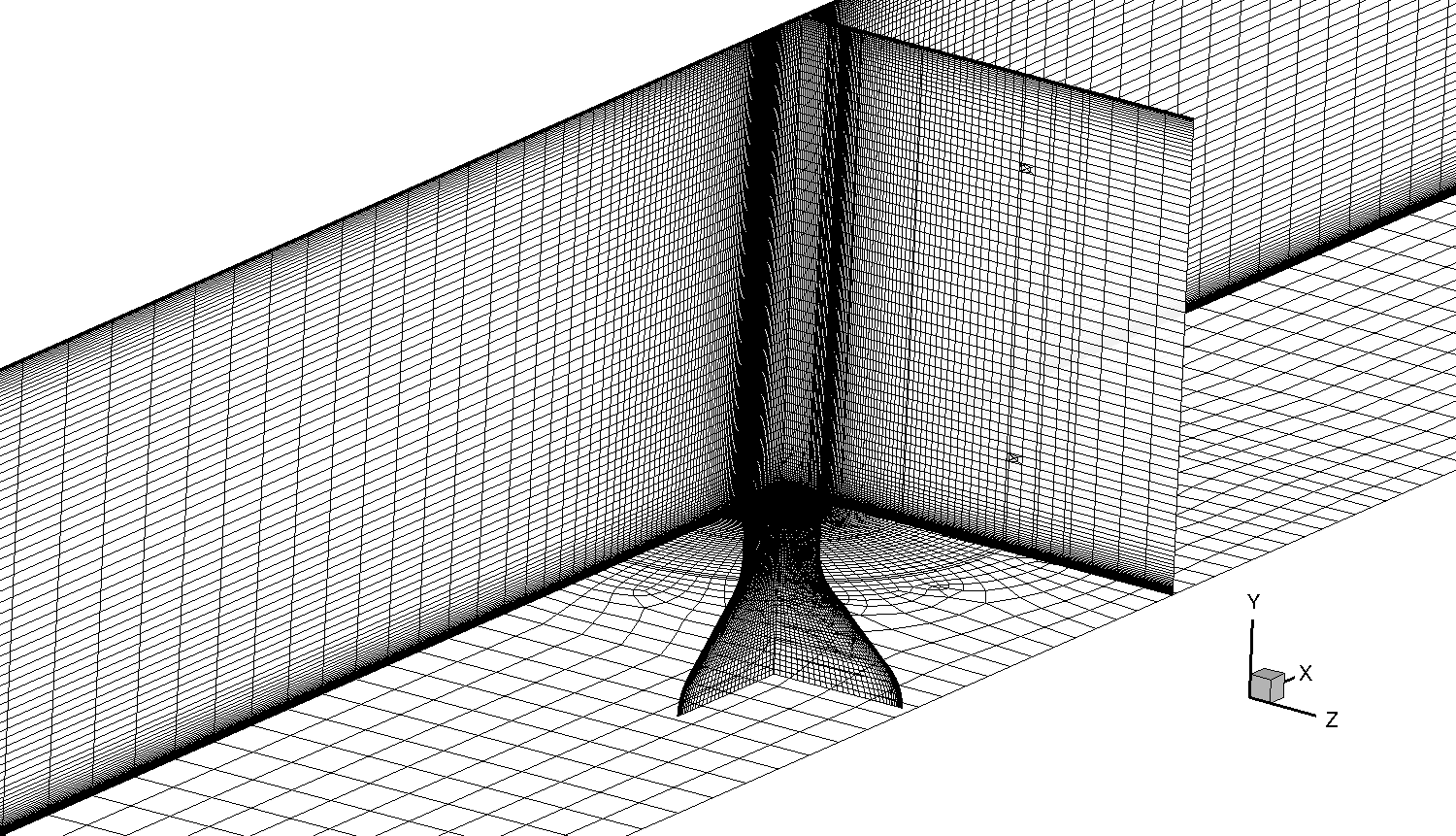}\label{fig:hyshot-3D-geom-c}}
\caption{The 3D combustion chamber mesh used for the HyShot II simulations, highlighting (a,c) the injection port and (b) the nozzle.}
 \label{fig:hyshot-3D-geom}
 \end{center}
\end{figure}

The grid is again designed such that the wall-normal spacing adjacent to the wall is $y_1^+ = 1$, and the total number of control volumes in this domain is $1.2$M. To increase numerical stability, a corner radius of $0.06$mm is used at the fuel nozzle orifice; see Figure \ref{fig:hyshot-3D-geom-a}. The domain is essentially structured except near the injection port; see Figure \ref{fig:hyshot-3D-geom-c}.

\subsubsection{Coupling the two domains} 
\begin{figure}[h!tb]
 \begin{center}
 \includegraphics[width=0.75\linewidth]{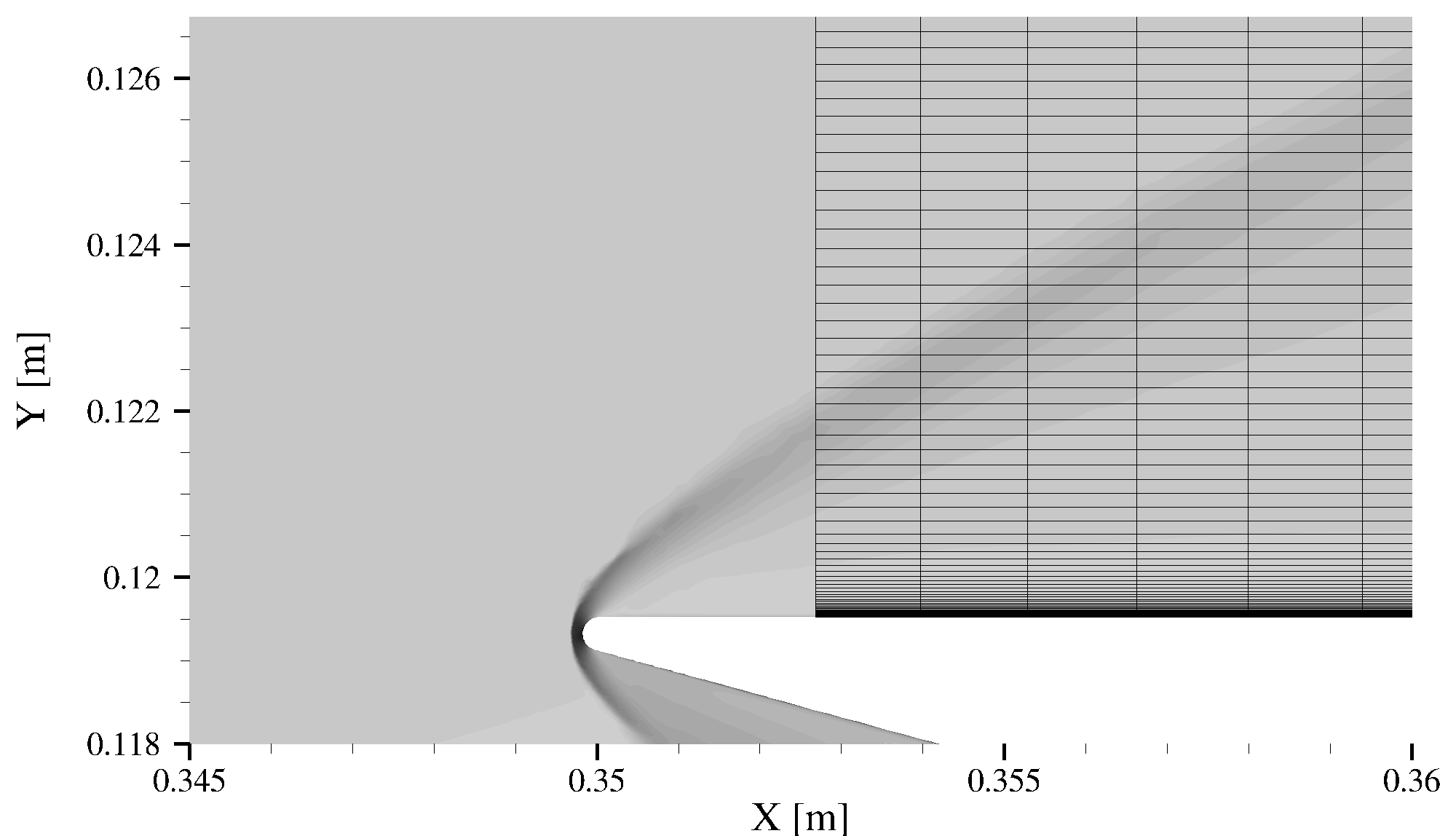}
 \caption{Density contours from the 2D domain overlaid by the 3D mesh (black lines) showing where the 3D inflow profile is extracted relative to the shock structures.}
 \label{fig:hyshot-mesh-overlap}
 \end{center}
\end{figure}

\noindent In the 3D domain the inflow condition for the oxidizer stream is taken from the 2D simulation. A wall-normal profile is extracted at $x=352.68$mm (absolute coordinates of the mesh), this profile is applied uniformly across the span of the combustion chamber inlet. The profile is extracted at a location where the oblique shock generated by the body wall leading edge is captured well above the boundary-layer resolving cells, see Figure \ref{fig:hyshot-mesh-overlap}.

\subsection{Reynolds-averaged Navier-Stokes}

\noindent Simulating the HyShot II scramjet with large-eddy simulation (LES) is tractable but costly~\cite{Chapius2013,larsson:14}. Wall-modeled LES uses approximately $100$M cells, and wall-resolved LES would require $100$B cells. It is not feasible to use LES for uncertainty quantification, which requires many simulations to estimate statistics of the solution behavior. Reynolds-averaged Navier-Stokes (RANS) costs less than LES by modeling only the time-averaged flow. Also, the mesh can be much coarser than LES, which further reduces the cost. The trade-off is that turbulence is now modeled; the accuracy of the result depends strongly on the accuracy of the turbulence model. RANS represents the most practical approach for performing UQ of the HyShot II system, and in general RANS is the most common approach used by hypersonic vehicle designers for simulating turbulence in their systems \cite{Yentsch2012,Karl2008}.

Pe\v{c}nik et al.~\cite{Pecnik2012} and Terrapon et al.~\cite{Terrapon2010} give comprehensive treatments of the discretization and numerical implementation of the corresponding physical models. Quantitative comparisons between the computational
predictions of the pressure within the HyShot II combustion chamber are reported in Pe\v{c}nik et al.~\cite{Pecnik2012}; these comparisons show remarkable agreement in terms of the compression ratio due to fuel injection and burning. Moreover, shock wave locations and strength within the chamber compare favorably. Overall, the simulations capture the dominating effects of the turbulent boundary layer interacting with the shocks. Furthermore, the results compare wall heat flux with experiments, again demonstrating satisfactory agreement. Pe\v{c}nik et al.~\cite{Pecnik2012} also present detailed grid resolution analysis, and they study of the effects of spanwise domain size and side walls; the solver and the results presented in~\cite{Pecnik2012} have been used to design the simulations presented in this paper.

\subsubsection{Physics modeling and closures}
\noindent Several RANS closures are needed in the multiphysics model. We next describe our specific choices.

\paragraph{Turbulence}
The turbulence model used to determine the Reynolds stresses is based on the eddy-viscosity hypothesis and the $k\!-\!\omega$ shear stress transport (SST) formulation~\cite{Menter93}. The SST model is one of the most frequently used RANS models in industrial applications. In many applications with compressible boundary layers and shock-turbulence interactions, SST gives reasonably accurate predictions \cite{Emory2013, oliver2007assessment}.
In the SST model, two transport equations are solved in addition to the RANS equations to describe the turbulent kinetic energy $k$ and the specific dissipation rate $\omega$ (having units of inverse time). The SST model blends the standard $k\!-\!\omega$ and $k\!-\!\epsilon$ models to accurately represent both the near wall regions and the response to high strain and pressure gradient.

Due to its popularity, several modified SST models have been proposed. We apply two specific limiter modifications~\cite{Pecnik2011}. The first is a limiter for the eddy-viscosity~\cite{menter1994},
\begin{equation}
  \mu_t = \frac{\rho a_1 k}{\max \left(a_1 \omega, \Omega F_2 \right)}, 
\end{equation}
where $\rho$ is the density, $a_1=0.31$ is a model constant, $k$ is the turbulence kinetic energy, $\omega$ is the specific dissipation, and $\Omega$ is the vorticity magnitude. $F_2$ is a function designed to be 1 for boundary-layer flows and 0 for shear layers. The second is a turbulence kinetic energy production limiter, which is meant to ensure realizable Reynolds stresses,
\begin{equation}
  \mathcal{P}_k = \min\left( \mathcal{P}_k, 20 \, C_\mu \, k \, \omega \right) \, ,
\end{equation}
where $\mathcal{P}_k$ is the turbulence kinetic energy production, $C_\mu=0.09$ is a model constant, and $k$ and $\omega$ as above.

\paragraph{Transition}
The laminar-turbulent boundary-layer transition process is modeled by manually specifying transition locations on the intake ramp, body wall, and cowl wall. These locations then affect the flow through inhibited production and destruction of turbulence kinetic energy in the boundary layers upstream of these points;
this mimics the handling of transition in the $\gamma\!-\!Re_{\theta t}$ transition model by \citet{Menter2004}. 

\paragraph{Combustion}

\noindent The combustion model is based on a flamelet/progress variable approach (FPVA), in which the chemistry is tabulated as a series of laminar flamelet solutions for a given set of boundary conditions and background pressure. The effect of turbulence on the flame is approximated by using a presumed beta probability density function  for the fuel/air mixture fraction. In this manner the chemical composition is mapped \textit{a priori} with respect to a small number of parameters used to search this table. The major assumption behind this approach is that chemistry is fast relative to the mixing time scales and can therefore be accurately represented by a small number of scalar quantities. This approach requires three additional transport equations for the mixture fraction, mixture fraction variance, 
and progress variable.
These are the values used with the table to provide species mass fractions and other properties that in turn influence the local temperature and pressure \cite{Pecnik2010,Terrapon2010,Pecnik2012}. The 20 reaction ${\rm H2\!-\!O2}$ mechanism of \citet{Hong2011} is used to generate the FPVA table. The FPVA framework is fairly common in subsonic application, as it has been developed based on a low Mach assumption. Formal descriptions of the development and extension to high speed flows are found in \citet{Terrapon2009,Saghafian2011}.

\paragraph{Boundary conditions}
In both the 2D and 3D domains, solid walls are modeled as isothermal with $T_w = 300$K due to the very short test time which
prevents the steel walls from heating up \cite{Gardner2004}.
The nozzle (combustor exit), bleed channel, and freestream (flow outside cowl side of vehicle) use a Neumann boundary condition.
At the fuel inlet, the stagnation pressure $P_{0,\rm H2}$ and the stagnation temperature $T_{0,\rm H2}$ are specified. The latter is always 300K, whereas the former is varied depending on the desired fuel equivalence ratio.

\subsubsection{Flow solver}
\label{sec:flowSolver}

\noindent The present calculations solve the steady, compressible RANS equations (five PDEs), the $k\!-\!\omega$ SST turbulence model (two PDEs), and the FPVA combustion model (three PDEs and a look-up table). The computations are carried out under the steady-state assumption although the unstart process is transient. This choice is motivated by the desire to detect the conditions that lead to unstart rather than model the entire unstart process. The flow solver \verb=Joe=, developed at Stanford's Center for Turbulence Research, is used to perform these simulations. The code performs parallel calculations on a collocated unstructured mesh using a finite volume formulation. The discretization is second order in space; gradients of flow quantities are calculated with least-squares. Comprehensive discussions of the numerical implementation as well as validation of \verb=Joe= for HyShot II simulations can be found in~\citet{Pecnik2010,Terrapon2010,Pecnik2012}.

\paragraph{Cost}
The average time for each 2D forebody simulation (Section \ref{sec:2d}) used in the uncertainty quantification study was approximately 8.1 hours on 48 CPUs. The average time for each 3D combustor simulation (Section \ref{sec:3d}) was approximately 76 hours on 120 CPUs, though there was significant variability in the number of iterations needed to meet the desired stopping criteria.

\section{Uncertainty sources and the quantity of interest}
\label{sec:sources}

\noindent The first step in uncertainty quantification is to identify and mathematically describe the system's uncertainties. These uncertainties must be propagated through the model to assess their impact on the quantity of interest. In this section we introduce and characterize both the sources of uncertainty and the quantity of interest in the HyShot II model. We first identify reasonable ranges for the input parameters. All available knowledge (e.g., from observations, theory, and expert opinion) is used to inform the input ranges. These ranges are sufficient to estimate (i) minimum and maximum values for the output quantity of interest and (ii) safe sets of operating conditions. To estimate a complete cumulative distribution function on the output quantity of interest, it is sufficient to consider a density function on the input space. We follow Jaynes' maximum entropy principle~\cite{jaynes1957information} to choose a uniform density on the space of inputs bounded by the ranges. Without any additional information (e.g., from actual measurements), the maximally ignorant or least informed assumption is that of a uniform density, which maximizes the (mathematical) entropy over all possible densities on the hyperrectangle defined by the ranges. The uniform density is a modeling choice, and the uncertainty quantification results depend on this choice. A complete study should analyze the results' sensitivity to the choice of density, but this is outside the scope of the present work. The parameters and associated ranges are summarized in Table \ref{tab:inflow_uq_ranges}.


\subsection{Inflow conditions, mean quantities}

\noindent The objective is to model conditions from the ground test experimental campaign carried out by researchers at DLR in the HEG. The conditions in the HEG represent the flight test from \citet{Boyce2003} along with boundary conditions and comprehensive measurements in the combustor. Stanford University and DLR collaborated to share experimental results, simulation results, and additional unpublished information. When this study was initialized, a total of thirteen experimental runs in the HEG shock tunnel at the relevant conditions (i.e., mimicking the high altitude conditions encountered during the flight test) had been performed, four without fueling and nine with fueling. For each experimental run, the stagnation pressure $P_0$, temperature $T_0$, and enthalpy $H_0$ of the nozzle supply region are measured (reported in Table \ref{tab:dlr-heg-shots}). The fueled runs additionally measured the hydrogen fuel plenum pressure $P_{0,{\rm H2}}$, from which the fuel/air equivalence ratio $\phi$ was estimated as
\begin{equation}
\label{eq:phi}
  \phi = \frac{8 \dot{m}_\textrm{H2}}{\dot{m}_\textrm{O2}} \,.
\end{equation}
The equivalence ratio is a natural quantity to consider given its direct chemical interpretation, but in this study we will consider the fuel plenum pressure as the control parameter since this was the directly measured quantity in the experiments.
For all runs the hydrogen plenum temperature was 300K.

\begin{table}[!t]
\begin{center}
\begin{tabular}{c l || c | c | c | c | c}
 & Shot & $P_0$ $\left[{\rm bar}\right]$ & $T_0$ $\left[{\rm K}\right]$ & $H_0$ $\left[{\rm MJ}/{\rm kg}\right]$ & $P_{H2}$ $\left[{\rm bar}\right]$ & $\phi$ \\ \hline
 \multirow{4}{*}{fuel-off}
  & 805 & 178.05 & 2742 & 3.25 \\
  & 807 & 181.19 & 2777 & 3.30 \\
  & 808 & 175.27 & 2716 & 3.21 \\
  & 814 & 176.64 & 2735 & 3.24 \\ \hline
  \multirow{9}{*}{fuel-on}
  & 804 & 172.96 & 2652 & 3.22 & 5.41 & 0.341 \\
  & 809 & 177.80 & 2753 & 3.27 & 5.27 & 0.329 \\
  & 810 & 176.66 & 2705 & 3.20 & 5.73 & 0.351 \\
  & 811 & 178.81 & 2726 & 3.23 & 4.68 & 0.286 \\
  & 812 & 179.84 & 2729 & 3.23 & 5.32 & 0.325 \\
  & 816 & 173.06 & 2769 & 3.39 & 3.71 & 0.266 \\
  & 817 & 176.99 & 2796 & 3.13 & 5.28 & 0.315 \\
  & 827 & 170.23 & 2701 & 3.19 & 5.09 & 0.324 \\
  & 828 & 187.43 & 2796 & 3.28 & 4.84 & 0.288
\end{tabular}
\end{center}
\caption{Experimental conditions from HyShot II ground tests performed in the HEG \cite{DLRReport}.}
\label{tab:dlr-heg-shots}
\end{table}

\begin{figure}[!t]
 \begin{center}
 \includegraphics[width=0.75\linewidth]{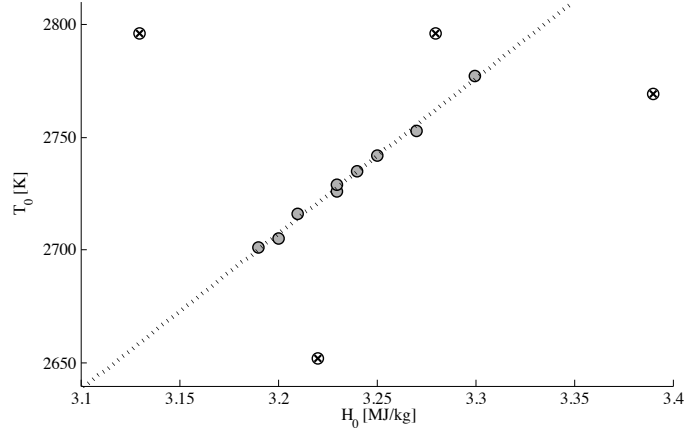}
 \caption{Least-squares fitting (dotted line) of $T_0\textrm{ vs. }H_0$ data from Table \ref{tab:dlr-heg-shots}, where the experiments marked with $\otimes$ (runs 804, 816, 817, and 828) have been excluded from the curve fit.}
 \label{fig:T0H0_relation}
 \end{center}
\end{figure}

Table~\ref{tab:dlr-heg-shots} reports the measured and estimated parameters from the 13 experimental runs. It is important to point out that the stagnation temperature $T_0$ and stagnation enthalpy $H_0$ are dependent parameters, being related by the specific heat $c_p$. When plotting these quantities in
Figure \ref{fig:T0H0_relation},
it is clear that 4 of the data points are contaminated by errors from some unknown source. The remaining data points (from 9 runs) fall closely on a straight line, which is consistent with the expectation of an approximately constant $c_p$ under these conditions. The least-squares regression model for these 9 data points is
\begin{equation}
\label{eq:T0H0-relation}
H_0 = \frac{T_0 - 508.1386}{6.8718\times10^{-4}} \,.
\end{equation}
We thus take the $P_0$ and $H_0$ values from Table \ref{tab:dlr-heg-shots}, and use \eqref{eq:T0H0-relation} to calculate the associated $T_0$. We transform these stagnation conditions into the static pressure, temperature, and velocity to use as inflow conditions in the simulation. This is done using the following ratios provided by DLR
\begin{subequations}
\label{eq:hyshot-ratios}
\begin{alignat}{1}
  \frac{P}{P_0}  &= 1.16 \times 10^{-4} \, , \\
  \frac{T}{T_0} &= 0.0978 \,  , \\
  \frac{U_{\rm mag}}{\sqrt{H_0}} &= 1.332 \, ,   \label{eq:hyshot-ratios-U}
 \end{alignat}
\end{subequations}
which are assumed valid for any run condition. Equation \eqref{eq:hyshot-ratios-U} gives the velocity magnitude, which together with the angle-of-attack $\alpha$ allows us to compute velocity components.

While the simulation takes the velocity components and the static pressure and temperature to specify the inflow condition, it is more consistent to characterize the uncertainty directly in the stagnation conditions ($P_0$ and $H_0$) that are measured in the experiments. The small number of measurement data available (the 13 runs reported in Table \ref{tab:dlr-heg-shots}) makes it difficult to characterize this uncertainty. Nevertheless, 13 samples in the table are used to compute sample means ($\widehat{P}_0$ and $\widehat{H}_0$) and standard deviations ($\widehat{\sigma}_{P0}$ and $\widehat{\sigma}_{H0}$). The uncertainty range is conservatively defined as $\widehat{P}_0 \pm3 \cdot \widehat{\sigma}_{P0}$ and similarly for $H_0$.

We also take the angle-of-attack $\alpha$ as uncertain, due to
both the physical mounting of the model in the tunnel and the flow-induced deviation from this angle during a run. The nominal angle-of-attack was $3.6^\circ$. Based on the expert opinion of an experienced experimentalist \citetext{M. Gamba, personal communication, July $5^{th}$, 2012}, we believe $\alpha$ can only be specified within $\pm1^\circ$. Additional evidence suggests this value is uncertain based on analysis performed during the construction of the HyShot II model \cite{Gardner}. The weight of the model and associated instrumentation caused static load deflection of the model up to $0.465^\circ$, a value which will most likely increase during dynamic loading. Also, the entire building shakes while HEG is operating, which further supports the choice model $\alpha$ as uncertain. One may reasonably assume that the model or test section is at an effective angle-of-attack relative to the flow (we ignore any possible uncertainty in the yaw angle of the vehicle). Thus the uncertainty range is specified as $3.6\pm1^\circ$. Table \ref{tab:inflow_uq_ranges} summarizes the nominal value and endpoints of the uncertainty range for these parameters.

\subsection{Inflow conditions, turbulence quantities}
\noindent The $k\!-\!\omega$ SST model requires specification of the turbulence kinetic energy $k$ and the specific dissipation-rate $\omega$ at the inflow. This is done by specifying a turbulence intensity $I$ and turbulence dissipation length scale $L_{t,\omega}$ (characteristic of the energy-carrying eddies). Thus,
\begin{subequations}
\label{eq:hyshot-turb-inflow}
\begin{alignat}{1}
 \label{eq:k_inlet}
 k = \frac{3}{2}(U_{\rm mag}\,I)^2 \, , \\
  \label{eq:omega_inlet}
 \omega = \frac{\sqrt{k}}{C_\mu^{1/4}L_{t,\omega}} \, ,
 \end{alignat}
\end{subequations}
where $U_{\rm mag}$ is the mean velocity magnitude and $C_\mu = 0.09$ \cite{Tu2008}.

There are no direct measurements the turbulence intensity and length scale---neither their nominal values nor their deviation from nominal. We therefore estimate both the nominal values and the uncertainty. Assumed constant stagnation pressure and pressure fluctuation data---specifically that $P'_{\rm rms}/P_0 = 10^{-5}$---are used to estimate the nominal turbulence intensity as $I=1\%$. Combining this with expert opinion \citetext{M. Gamba, personal communication, July $5^{th}$, 2012}, we arrive at a range of $I = [0.1,1.9]$.

The estimates of turbulence length scale are informed by expert opinion \citetext{M. Gamba and J. Larsson, personal communication, July $5^{th}$, 2012}. The largest turbulence length scales generated at the throat preceding the nozzle in the tunnel
are roughly half the size of the throat. Furthermore, 3D isotropic eddies should grow as the inverse of the density in the isentropic expansion in the nozzle,
\begin{equation}
\label{eq:length_ratio}
 \frac{L_t}{L_{t,{\rm throat}}} \approx \left(\frac{P_0}{P}\frac{T}{T_0}\right)^{\frac{1}{3}}.
\end{equation}
The flow in the nozzle should approximately satisfy the 1D variable-area flow relation
\begin{equation}
 \left(\frac{A}{A^*}\right)^2 = \frac{1}{M^2}\left[\frac{2}{\gamma+1}\left(1+\frac{\gamma-1}{2}M^2\right)\right]^{(\gamma+1)/(\gamma-1)} \, ,
\end{equation}
where $\gamma$ is the ratio of specific heats. To achieve Mach $7.4$ (the nominal HEG flow condition) the nozzle must have an area ratio of $\frac{A}{A^*}=133$. The diameter of the test section at the HyShot II model is approximately $610$mm, thus the throat diameter is $53$mm. Substituting the pressure and temperature ratios from  \eqref{eq:hyshot-ratios} into \eqref{eq:length_ratio}, we find a length scale ratio of $9.43$, which leads to a nominal length scale $L_{t,\omega 0} \approx 245$mm. A conservative uncertainty range is specified relative to this nominal value where $L_{t,\omega} = L_{t,\omega 0}\pm50\%$. Table \ref{tab:inflow_uq_ranges} summarizes the nominal value and endpoints of the uncertainty range for these parameters.

\subsection{Transition locations}

\noindent There are three transition locations within the HyShot II: along the intake ramp/forebody (location labeled \xt{r}), along the body wall (\xt{b}), and along the cowl wall (\xt{c}). Along the forebody \xt{r} is determined from heat flux measurements at the vehicle surface. The resolution of thermocouples is very coarse, however, with only four thermocouples spanning $115$mm \cite{DLRReport}. The body and cowl transition locations are specified through expert opinion (based on the experience of the HEG facility operators) and there are no transition measurements available on these walls. Based on this information, the three transition locations are assumed uncertain---though are still treated as fixed locations (ignoring uncertainty related to the assumption that these locations are stationary).

In this work we consider uncertainty in \xt{r} and \xt{c}, both of which reside in the 2D computational domain. Note that \xt{r} is defined relative to the vehicle leading edge and \xt{c} is defined relative to the cowl leading edge; these locations are shown in Figure \ref{fig:hyshot-2D-geom}. The 2D simulations are significantly cheaper than their 3D counterpart, so the impact of considering additional uncertainties is relatively small. The 3D computations are more expensive due to both increased domain size and the additional modeling closures (e.g., combustion, mixing). Thus uncertainty in \xt{b} is ignored as a practical consideration to reduce computational cost of the uncertainty quantification.

DLR uses a simple transition model to estimate the location \xt{r}=145mm because of the limited resolution provided by the thermocouples mounted on the forebody walls; this estimate is consistent with experimental correlations that indicate transition occurs where
\begin{equation}
  \label{eq:trans_crit}
  \frac{Re_\theta}{M_e} \approx 200 \, ,
\end{equation}
where $Re_\theta$ is the Reynolds number based on momentum thickness and $M_e$ is the Mach number at the edge of the boundary-layer. As opposed to directly assuming that the magnitude of \xt{r} is uncertain, we instead assume that the critical value of $Re_\theta/M_e=200$ from transition theory is uncertain.

It is well known (see, e.g.,~\cite{Menter2004}) that the transition Reynolds number is sensitive to several factors including the free-stream turbulence, the wall roughness, etc. It is therefore difficult to specify a precise threshold value. In our previous uncertainty analysis~\cite{Pecnik2011}, we consider the effect of compressibility on the growth of turbulent spots in the transition region by introducing uncertainty in the correlation specification in \eqref{eq:trans_crit}. Here we follow the same strategy and define that transition is equally likely to occur
at the point where
\begin{equation}
\label{eq:trans_critical}
 \frac{Re_\theta}{M_e}=200\,(1 \pm \varphi),
\end{equation}
where $\varphi=0.2$. The ratio $Re_\theta/M_e$ is linear in the momentum thickness $\theta$, which grows as $\theta(x) \propto \sqrt{x}$ in the laminar boundary layer prior to transition. Therefore, we can write
\begin{equation}
 \frac{Re_\theta}{M_e} = c \, \sqrt{x} \,,
\end{equation}
for some constant $c$. Linearizing about the nominal value \xt{r0} gives
\begin{alignat}{1}
 \frac{Re_\theta}{M_e} &= c\, \sqrt{x_{t,{\rm r}0}} + \frac{c}{2 \sqrt{{x}_{t,{\rm r}0}}}(x-{x}_{t,{\rm r}0}) \, , \nonumber \\
 &= c \, \sqrt{{x}_{t,{\rm r}0}} \left(1 \pm \frac{x-{x}_{t,{\rm r}0}}{2{x}_{t,{\rm r}0}}\right) \, .
\end{alignat}
This is in the same form as \eqref{eq:trans_critical}, thus the uncertainty $\varphi$ is
\begin{equation}
 \varphi = \pm \frac{x-{x}_{t,{\rm r}0}}{2{x}_{t,{\rm r}0}} \, .
\end{equation}
Rearranging we can represent the transition location uncertainty as a function of $\varphi$
\begin{equation}
\label{eq:transx}
 {x}_{t,{\rm r}} = {x}_{t,{\rm r}0}(1 \pm 2\varphi) \, ,
\end{equation}
thus the uncertainty range is ${x}_{t,{\rm r}} = {x}_{t,{\rm r}0}\cdot [0.6,1.4]$ where ${x}_{t,{\rm r}0} =145$mm.

The same procedure is used to determine the uncertainty in \xt{c}. The nominal transition location is specified as $50$mm downstream of the cowl leading edge, so the range is defined as ${x}_{t,{\rm c}} = {x}_{t,{\rm c}0}\cdot [0.6,1.4]$ where ${x}_{t,{\rm c}0} = 50$mm. The uncertainty ranges and nominal values for these parameters are summarized in Table \ref{tab:inflow_uq_ranges}.

\begin{table}[h!tb]
 \begin{center}
 \begin{tabular}{c | c c c | c}
 Parameter    & Min  & Nominal & Max & Units   \\ \hline 
   Stagnation Pressure      & 16.448   & 17.730  & 19.012  & MPa     \\
   Stagnation Enthalpy        & 3.0551   & 3.2415  & 3.4280  & $\mbox{MJ}/\mbox{kg}$ \\
   Angle of Attack     & 2.6      & 3.6     & 4.6     & deg. \\
   Turbulence Intensity          & 0.001    & 0.01    & 0.019   & $\cdot$ \\
   Turbulence Length Scale          & 0.1325   & 0.245   & 0.3575  & m       \\
   Ramp Transition Location & 0.087    & 0.145   & 0.203   & m       \\
   Cowl Transition Location & 0.030    & 0.050   & 0.070   & m       
\end{tabular}
\end{center}
 \caption{Summary of parameters and ranges used in the uncertainty quantification of the HyShot II scramjet within the 2D intake ramp simulation.}
 \label{tab:inflow_uq_ranges}
\end{table}

\subsection{Quantity of interest}
\label{sec:rans:QoI}

\noindent Numerical simulations where the shock system has moved all the way to the upstream boundary
(i.e., when the flow is about to unstart)
are computationally challenging. For this reason, we seek a proxy indicator for the unstart process.
Several possible unstart indicators have been proposed \cite{Emory2011}. Presently, we use a function of the combustion chamber exit pressure. The motivation for this choice is related to the physics of compressible flow with heat addition (Rayleigh's flow): a supersonic flow decelerates towards sonic conditions when heat is released through a combustion process. The corresponding increase in pressure achieves its maximum just before choking conditions are reached.


The unstart proxy is the normalized integral of pressure over the last $10$mm of the combustor
\begin{equation}
  \label{eq:unstartProxy}
  \frac{1}{\mbox{vol}\left(\boldsymbol{V}\right)} \iiint \limits_{\boldsymbol{V}} P \; dx\; dy\; dz \, ,
\end{equation}
where $\boldsymbol{V}$ is the combustor region in $0.64 \leq x \leq 0.65$m. The integrated exit pressure \eqref{eq:unstartProxy} is a critical measure of the scramjet performance, and we treat it as the quantity of interest for uncertainty quantification. 




\section{Uncertainty quantification methodology}
\label{sec:uq}

\noindent In this section, we use the notation $f(\vx)$ to represent the quantity of interest (integrated exit pressure defined in \eqref{eq:unstartProxy}) as a function of the input parameters (Table \ref{tab:inflow_uq_ranges}), where we treat the parameters as normalized to the hypercube $[-1,1]^m$ ($m=7$ in the HyShot II model). In other words, for a point $\vx\in[-1,1]^m$, the function $f(\vx)$ first shifts and linearly scales $\vx$ to the application-specific ranges and then computes the quantity of interest from the given simulation with the proper inputs. This notation simplifies the methods' description. 

The uncertainty quantification problem for the HyShot II model involves three distinct computations: (i) estimating lower and upper bounds for $f(\vx)$, (ii) identifying sets of inputs that produce safe operation defined as $f(\vx)$ remaining below a given threshold, and (iii) estimating a cumulative distribution function for $f(\vx)$ given a probability density function on $\vx$. The last computation is the typical \emph{forward} uncertainty quantification problem~\cite{smith2013uncertainty}. These computations are challenging for several reasons.
\begin{enumerate}
\item The exit pressure $f$ is a functional of the pressure field computed from the complex, nonlinear, multiphysics simulation, so we have no prior knowledge of exploitable structure like linearity or convexity.
\item Each evaluation of $f(\vx)$ requires an expensive simulation; each run takes approximately 9500 CPU-hours on available computing resources.
\item Due to the multiphysics and domain coupling, we cannot evaluate gradients or Hessians of $f(\vx)$ with respect to $\vx$. 
\item Evaluations of $f(\vx)$ contain nonnegligible numerical error due to (i) the fixed point iteration scheme that solves the compressible flow model and (ii) the fixed mesh constructed to capture shocks at the nominal inputs.
\end{enumerate}
Under these conditions, exhaustively sampling the seven-dimensional input space is infeasible for our three desired computations. Reducing the simulation's cost with parametric reduced-order models~\cite{Benner13asurvey} is infeasible for the nonlinear, multiphysics simulation code. Multivariate approximation schemes for response surfaces---e.g., splines or radial basis functions~\cite{wendland2004scattered}---are precarious since we cannot densely sample the input space. Adaptive methods for response surface construction use heuristics (e.g., \emph{expected improvement}) to increase efficiency~\cite{jones2001taxonomy}, but their efficacy relies on a sufficiently accurate initial approximation, which is not guaranteed in high dimensions with limited samples. Sophisticated methods for response surfaces in high dimensions---e.g., sparse grids~\cite{Bungartz2004}, low-rank approximations~\cite{Doostan2013}, or sparse approximation schemes~\cite{Doostan2011}---are appropriate when $f(\vx)$ admits particular structure; we have no prior indication that the integrated exit pressure as a function of the HyShot II simulation inputs admits such structure. Global sensitivity metrics (e.g., Sobol indices) can identify the least important parameters, which we might fix at nominal values to reduce the dimension for the computations~\cite{saltelli2008global}. But such metrics require approximating seven-dimensional integrals; insufficient sampling can yield incorrect results.

\subsection{Active subspaces and summary plots}
\label{sec:as}

\noindent Our approach discovers and exploits an \emph{active subspace} in $f(\vx)$~\cite{constantine2013active}. The active subspace is defined by a set of important directions in the input space, which generalizes the idea of identifying the most important coordinates with sensitivity analysis. The important directions come from the first few eigenvectors of the $m\times m$ symmetric, positive semidefinite matrix
\begin{equation}
\label{eq:C}
\mC \;=\; \int \nabla f(\vx)\,\nabla f(\vx)\,\rho(\vx)\,d\vx \;=\; \mW\Lambda\mW^T,
\end{equation}
where $\nabla f$ is the gradient of $f$ with respect to $\vx$, $\rho$ is a probability density function on the input space, $\mW$ is the orthogonal matrix of eigenvectors, and $\Lambda$ is the diagonal matrix of nonnegative eigenvalues. The $i$th eigenvalue $\lambda_i$ satisfies
\begin{equation}
\lambda_i \;=\; \int (\vw_i^T\nabla f(\vx))^2\,\rho(\vx)\,d\vx,\quad i=1,\dots,m.
\end{equation}
In words, the eigenvalues measure how $f$ changes, on average, in response to small perturbations in $\vx$ along the corresponding eigenvectors. For $n<m$, if there is a large gap between the $n$th and $n+1$th eigenvalues, then a reasonable approximation for $f(\vx)$ is 
\begin{equation}
\label{eq:lowd}
f(\vx) \;\approx\; g(\mW_1^T\vx),
\end{equation}
where $\mW_1$ contains the first $n$ columns of $\mW$, and $g$ is a map from $\mathbb{R}^n$ to $\mathbb{R}$. The model defined by \eqref{eq:lowd} is a low-dimensional approximation that can be exploited to enable otherwise infeasible parameter studies~\cite{constantine2013active}. If the gradient $\nabla f(\vx)$ is available as a subroutine, then one can estimate $\mC$ and its eigenpairs in \eqref{eq:C} with Monte Carlo~\cite{constantine2014computing}. 

Our simulation code does not have a subroutine to compute the gradient. We do not trust finite difference approximations because of fixed point solver oscillations at nearby parameter values. To find the active subspace, we turn to methods for \emph{sufficient dimension reduction} (SDR) developed in the context of regression modeling~\cite{cook2009regression}. The simplest of these methods is based on a least-squares-fit linear approximation of $f(\vx)$. Assume that $\rho:\mathbb{R}^m\rightarrow\mathbb{R}_+$ is a given probability density function on the inputs $\vx$; as discussed in Section \ref{sec:sources}, we follow the maximum entropy principle to choose a uniform density on the hypercube $[-1,1]^m$ in the absence of information to the contrary. Consider the procedure outlined in Algorithm \ref{alg:as}. 

\begin{algorithm}{Active subspace estimation with least-squares linear approximation}{\label{alg:as}}
\begin{enumerate}
\itemsep0em
\item Draw $M$ samples $\{\vx_j\}$ independently according to the density $\rho(\vx)$.
\item For each $\vx_j$, compute $f_j = f(\vx_j)$.
\item Determine the parameters $\hat{u}_0,\dots,\hat{u}_m$ of the linear approximation
\begin{equation}
\label{eq:linapprox}
f(\vx) \;\approx\; \hat{u}_0 + \hat{u}_1x_1 + \cdots + \hat{u}_mx_m
\end{equation}
with least-squares
\begin{equation}
\label{eq:lscomp}
\hat{\vu}\;=\;\argmin{\vu}\; \frac{1}{2}\,\|\mA\vu - \vf\|_2^2,
\end{equation}
where 
\begin{equation}
\mA = \bmat{1 & \vx_1^T \\ \vdots & \vdots \\ 1 & \vx_M^T},\quad
\vf = \bmat{f_1\\ \vdots \\ f_M},\quad
\hat{\vu} = \bmat{\hat{u}_0 \\ \vdots \\ \hat{u}_m}.
\end{equation}
\item Compute the vector $\vw$ as
\begin{equation}
\label{eq:w}
\vw \;=\; \frac{\hat{\vu}'}{\|\hat{\vu}'\|},
\end{equation}
where $\hat{\vu}'=[\hat{u}_1,\dots,\hat{u}_m]^T$ is the last $m$ coefficients (i.e., the gradient) of the linear approximation \eqref{eq:linapprox}.
\end{enumerate}
\end{algorithm}

We define the active subspace with the unit vector $\vw$ from \eqref{eq:w}, which identifies one important direction in the input space. The number $M$ of samples must be $\mathcal{O}(m)$--- large enough to fit the linear approximation \eqref{eq:linapprox} with least-squares---and the sample points $\{\vx_j\}$ must produce a full rank least-squares problem. Some recent work studies the scaling of $M$ with the dimension of the polynomial space for optimal approximation of the coefficients~\cite{migliorati2014}; for a linear approximation \eqref{eq:linapprox}, the polynomial space dimension is $m+1$. 

Several remarks are in order. We do not use the linear approximation as a predictive response surface. The coefficients $\hat{u}_1,\dots,\hat{u}_m$ of the linear approximation identify a direction in the input space---hence the normalization in \eqref{eq:w}. The vector $\vw$ may reveal low-dimensional structure in $f(\vx)$; our goal is to exploit this structure, if present, to enable otherwise infeasible uncertainty quantification with the expensive computer simulation. 

In the SDR literature, Li and Duan~\cite[Theorem 2.1]{Li1989} show that under certain conditions, the least-squares-based procedure in Algorithm \ref{alg:as} identifies the dimension reduction space for a regression model. However, using the uniform density on the hypercube violates the assumptions of their theory. Additionally, treating the deterministic $f(\vx)$ with regression tools poses a challenge for interpretation; the standard regression model assumes random noise in the $f_j$'s that is not present in our deterministic simulation. Nevertheless, we can apply the procedure as outlined and validate that $\vw$ identifies the low-dimensional structure with a \emph{summary plot}~\cite{cook2009regression}.

The summary plot is a scatter plot of $\vw^T\vx_j$ versus $f_j$. The plot may reveal a nearly one-to-one mapping between the linear combination $\vw^T\vx$ and the output $f(\vx)$. In other words, the output may be well-represented by a scalar-valued function of the linear combination of $\vx$. Any perceived departure from the one-to-one mapping comes from one of two sources: (i) $f(\vx)$ may vary significantly as $\vx$ moves orthogonally to $\vw$ or (ii) the number $M$ of samples may be too small to accurately estimate $\vw$. In practice, one can determine whether the perceived departure is due to (i) or (ii) by running more simulations (i.e., increasing $M$). 

If the $f_j$'s depart very little from a one-to-one mapping, then we can reasonably approximate
\begin{equation}
\label{eq:1d}
f(\vx) \;\approx\; g(\vw^T\vx),\quad \vx\in[-1,1]^m,
\end{equation}
where $g:\mathbb{R}\rightarrow\mathbb{R}$ is a function constructed from the pairs $(\vw^T\vx_j,f_j)$. This construction is similar to a single-index regression model for the exit pressure with a quadratic polynomial link function~\cite[Chapter 11]{hastie2009elements}. Several constructions for $g$ are possible---e.g., polynomial approximation or radial basis approximation---but we prefer not to exactly interpolate the $f_j$'s due to finite sampling noise in the computed $\vw$ from \eqref{eq:w}. Note the similarity between \eqref{eq:1d} and \eqref{eq:lowd}; both approximate $f(\vx)$ with a function of less than $m$ linear combinations of $\vx$. 

The summary plot is a subjective tool, in general. The scientist perceives the $f_j$'s departure from a supposed univariate function of $\vw^T\vx$, and she must decide if this departure is small enough to justify the model \eqref{eq:1d}. One might use metrics such as the residual norm to compare different $g$'s (e.g., different degrees of polynomial approximation), but such metrics can be deceiving without the aid of visualization. We use a bootstrap method to gain some confidence in the perceived structure.

\subsection{Bootstrap for $\vw$}
\label{sec:bootstrap}

\noindent The elements of $\vw$ depend on the samples $f_j$ used to fit the linear approximation \eqref{eq:linapprox}. With $M$ random samples, it is natural to ask if the gradient of the linear approximation has been sufficiently resolved to reveal the true one-dimensional active subspace, if one is present. But a limited budget of function evaluations prohibits us from checking convergence of $\vw$ as $M$ increases. We use a bootstrap technique~\cite{Efron1993} to estimate the variability in the computed components of $\vw$ from the set of $M$ samples. 

\begin{algorithm}{Bootstrap for estimating variability in $\vw$}{\label{alg:bs}}
\begin{enumerate}
\itemsep0em
\item Choose the number $N$ of \emph{bootstrap replicates}.
\item For $k$ from 1 to $N$, do the following.
\begin{enumerate}
\itemsep0em
\item Let
\begin{equation}
\vpi^k \;=\; \bmat{
\pi_1^k & \cdots & \pi_M^k
}
\end{equation}
be $M$-vectors of integers between $1$ and $M$ sampled with replacement. 
\item Let the $M$-vector $\vf_k$ and $M\times(m+1)$ matrix $\mA_k$ be 
\begin{equation}
\vf_k \;=\; \bmat{f_{\pi_1^k} \\ \vdots \\ f_{\pi_M^k}},
\qquad
\mA_k \;=\; 
\bmat{1 & \vx_{\pi_1^k}^T \\ 
\vdots & \vdots  \\
1 & \vx_{\pi_M^k}^T}.
\end{equation}
\item Compute $\vw_k$ for each $\mA_k$ and $\vf_k$ as in \eqref{eq:lscomp} and \eqref{eq:w}. 
\end{enumerate}
\end{enumerate}
\end{algorithm}
The procedure outlined in Algorithm \ref{alg:bs} samples, with replacement, the rows of $\mA$ and elements of $\vf$ in \eqref{eq:lscomp} to create a collection of bootstrap replicates $\{\vw_k\}$. We can use these vectors to study the variability in the estimate $\vw$ from \eqref{eq:w}---e.g., with histograms that estimate the marginal bootstrap densities of $\vw$'s components. Sharp peaks and narrow supports in the histograms provide confidence in $\vw$. Wide supports may indicate a poor estimate of $\vw$, or it may reflect insufficient sampling of the $f_j$'s used to fit the linear approximation.

\section{Discovering and exploiting a one-dimensional active subspace in the HyShot II quantity of interest}
\label{sec:ashyshot}

\noindent We apply Algorithm \ref{alg:as} to the exit pressure \eqref{eq:unstartProxy} as a function of the $m=7$ input parameters in Table \ref{tab:inflow_uq_ranges}.
We repeat the study for two values of the fuel plenum pressure $P_{0,\rm H2}$, specifically
$P_{0,\rm H2}=4.8$ bar (corresponding to a fuel/air equivalence ratio of 0.30 at nominal inflow conditions)
and
$P_{0,\rm H2}=5.6$ bar (corresponding to a fuel/air equivalence ratio of 0.35 at nominal inflow conditions).
The first is the nominal operating condition. The second is chosen to be close to the critical equivalence ratio of 0.38-to-0.39. Recent studies found a change in the flow regime from an \emph{as-designed} fully supersonic combustion ($\phi < \sim 0.39$) to a state with a shock-train in the downstream part of the combustor (for $\phi > \sim 0.39$)~\cite{laurence:14,larsson:14}.
From an engineering point-of-view, it is interesting to see whether the effects of the inflow uncertainties change as one approaches this regime boundary. For $P_{0,\rm H2}=4.8$ bar, we use $M=50$ samples to estimate $\vw$, which is roughly $m^2$.
For $P_{0,\rm H2}=5.6$ bar, we use $M=14$ samples, which is $2m$. In both cases, the cost of searching for the one-dimensional active subspace is much less than the cost of estimating the integral quantities needed for global sensitivity analysis. 

\begin{figure}[htb]
\centering
\subfigure[$P_{0,\rm H2}=4.8$ bar]{\label{fig:as1}
\includegraphics[width=0.46\textwidth]{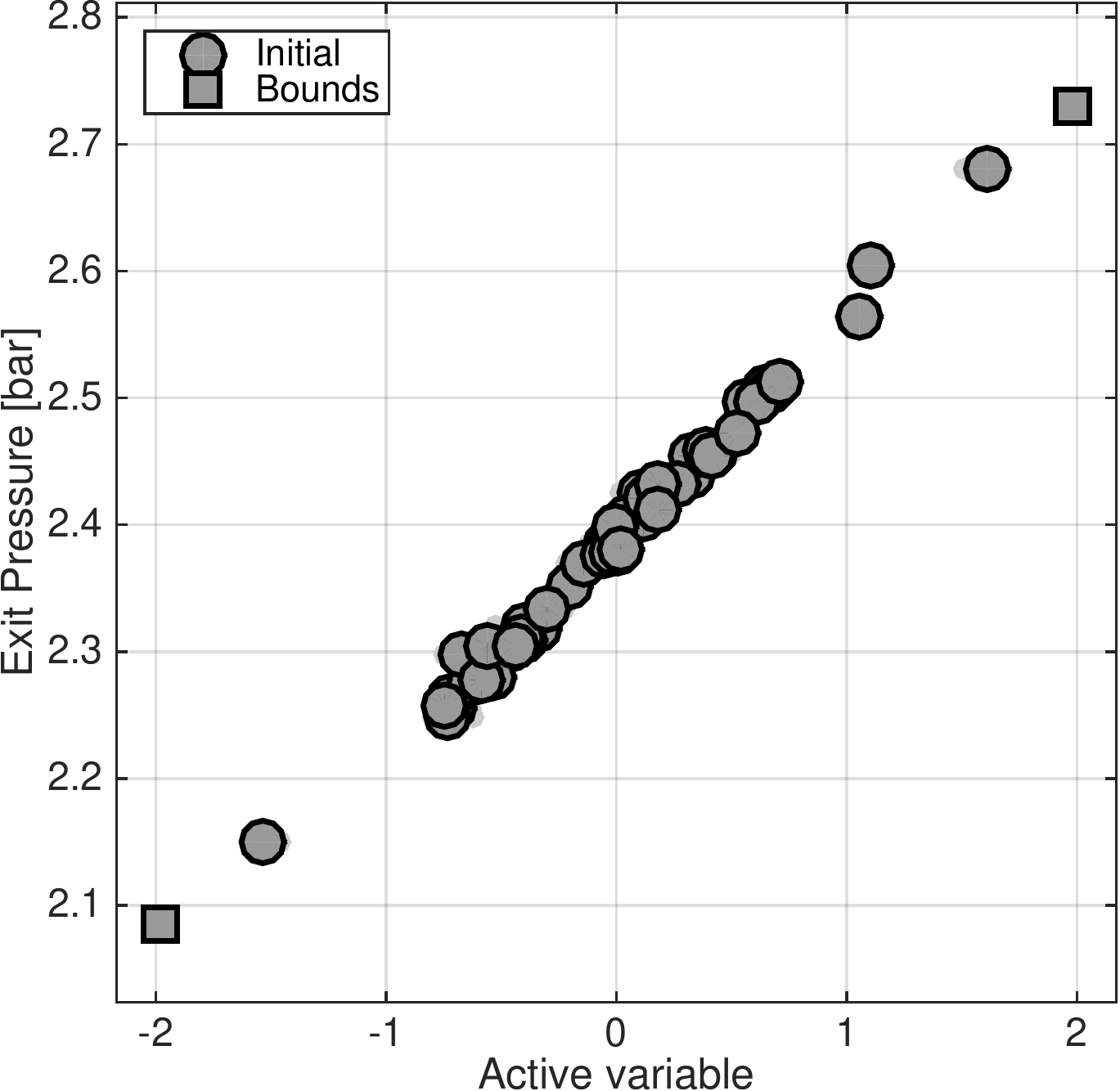}
}
\subfigure[$P_{0,\rm H2}=5.6$ bar]{\label{fig:as2}
\includegraphics[width=0.46\textwidth]{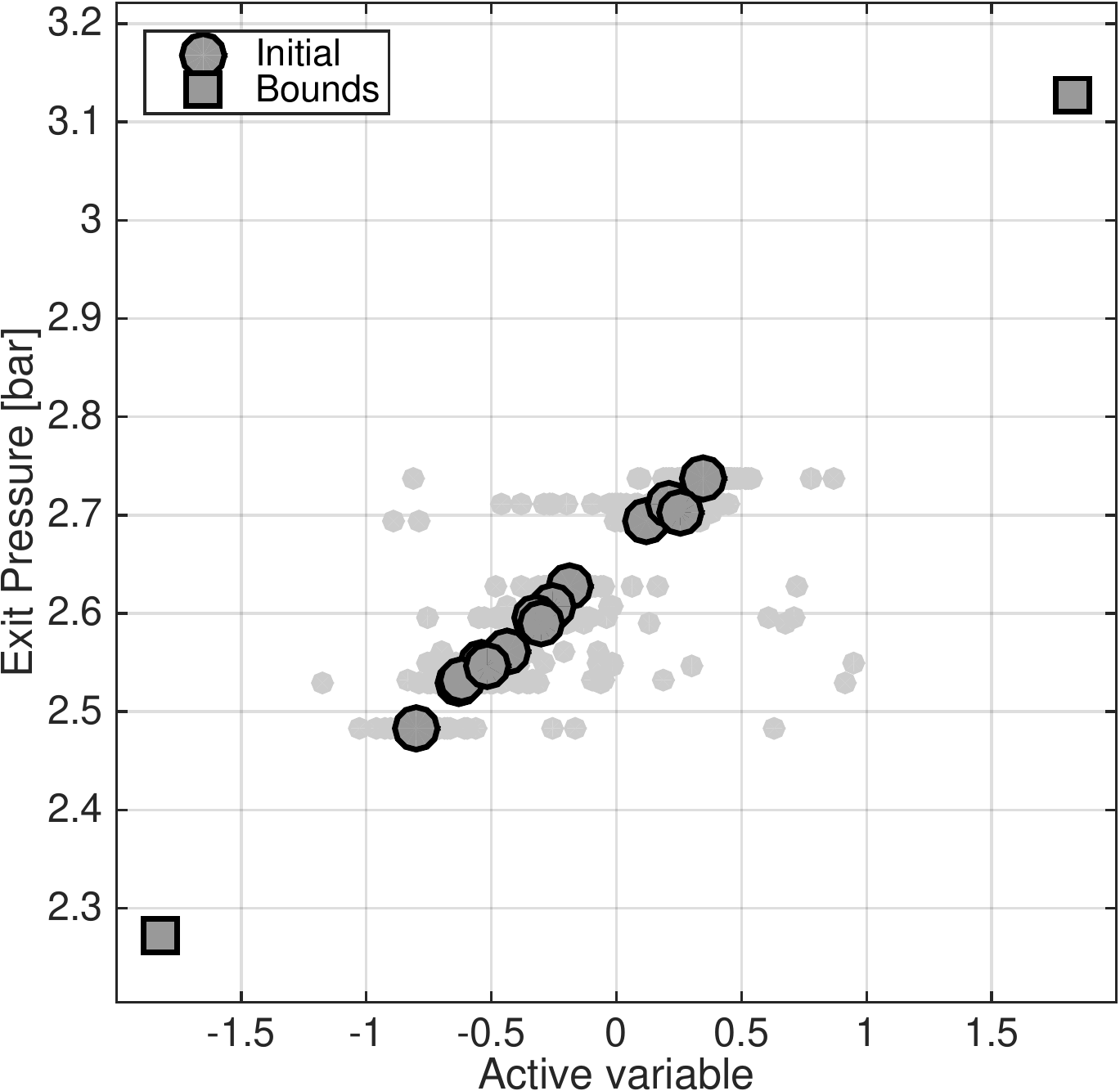}
}
\caption{The black circles are exit pressures computed from a set of HyShot II simulations and plotted against the corresponding active variable values $\vw^T\vx_j$. The left subfigure shows $M=50$ runs with fuel plenum pressure $P_{0,\rm H2}=4.8$ bar. The right subfigure shows $M=14$ runs with $P_{0,\rm H2}=5.6$ bar. A clear univariate relationship exists for both cases that we exploit to quantify uncertainty. The overlapping gray circles distributed horizontally correspond to exit pressures plotted against 100 bootstrap replicates of the active variable. There is hardly any variability in the bootstrap estimates for $P_{0,\rm{H2}}=4.8$. The horizontal variability for $P_{0,\rm H2}=5.6$ bar is a result of the relatively low oversampling ($M=2m$) when fitting the linear approximation \eqref{eq:linapprox}. The gray squares show the exit pressures at the boundaries of the domain where the perceived trend suggests we find the upper and lower bounds of the exit pressure. The univariate relationship is validated by these two additional simulations for each $P_{0,\rm H2}$.}
\label{fig:as}
\end{figure}

Figure \ref{fig:as} shows the summary plot for our application. The horizontal axis is the \emph{active variable} $\vw^T\vx$, and the vertical axis is the exit pressure. The black circles are the samples $f_j$ plotted against the corresponding value of the active variable $\vw^T\vx_j$. In our judgment, the circles indicate that we can approximate the exit pressure with sufficient accuracy using a univariate function of $\vw^T\vx$ as in \eqref{eq:1d} for both
$P_{0,\rm H2}=4.8$ bar and $P_{0,\rm H2}=5.6$ bar. We exploit this low-dimensional structure for the three desired computations needed for uncertainty quantification. 
\subsection{Sensitivity analysis}
\label{sec:sensitivity}

\noindent The components of $\vw$ measure the global sensitivity of the exit pressure to each of the seven parameters; see~\cite[Section 1.2.5]{saltelli2008global}. These values and their corresponding input parameters are shown in Table \ref{tab:dir}. They suggest that four of the seven parameters contribute the most to the one-dimensional active subspace. The angle of attack, the stagnation conditions (temperature and pressure), and the turbulence intensity dominate the transition location and the turbulence length scale. The largest contributors to changes in the exit pressure are the inputs related to the intensity of the bow shock at the nose of the vehicle, which influences the thermodynamic post shock conditions and the flow entering the Hyshot combustor. The state of the boundary layers (as determined by different transition location) plays a secondary role. It is more difficult to assess why the turbulence intensity plays such an important role. Caution is needed before interpreting this result, since RANS eddy-viscosity models (like the SST) are known to behave erratically across strong shock waves---like those present here \cite{sinha:03}; a more detailed investigation is outside the scope of the present paper. 

There is an important distinction between the plenum pressures $P_{0,\rm H2}=4.8$ bar and 5.6 bar. The leading factors in determining changes in the exit pressure do not change, but the relative importance of stagnation pressure and enthalpy does.
Large-eddy simulations of the HyShot {II} combustor~\cite{larsson:14} have shown that the combustion is close to complete, and thus the total heat release is essentially set by the fuel mass flux and thus independent of variations in the inflow conditions.
Those same large-eddy simulations also showed that the heat losses through the walls and the momentum loss due to friction are very large in this combustor.
Therefore, the most plausible explanation for the sensitivity of the exit pressure to the inflow stagnation pressure and enthalpy is that these variations must alter the friction and/or wall heat losses.

\begin{table}[htb]
\begin{center}
\begin{tabular}{cccc}
Index & $P_{0,\rm H2}=4.8$ bar & $P_{0,\rm H2}=5.6$ bar & Parameter\\
\hline
1 & 0.6506 & 0.7066 & Angle of Attack\\
2 & 0.5565 & 0.5008 & Turbulence Intensity\\
3 & -0.0002 & 0.0289 & Turbulence Length Scale\\
4 & 0.3685 & 0.2051 & Stagnation Pressure\\
5 & -0.3566 & -0.4490 & Stagnation Enthalpy\\
6 & -0.0196 & -0.0591 & Cowl Transition Location\\
7 & 0.0607 & 0.0432 & Ramp Transition Location
\end{tabular}
\caption{The components of the vector $\vw$ that defines the one-dimensional active subspace---computed independently for each $P_{0,\rm H2}$. Each component corresponds to one of the parameters in the HyShot II simulation as described in Section \ref{sec:sources}. }
\label{tab:dir}
\end{center}
\end{table}

\subsection{Bootstrap results}

\begin{figure}[htb]
\centering
\subfigure[Angle of Attack]{
\includegraphics[width=0.22\textwidth]{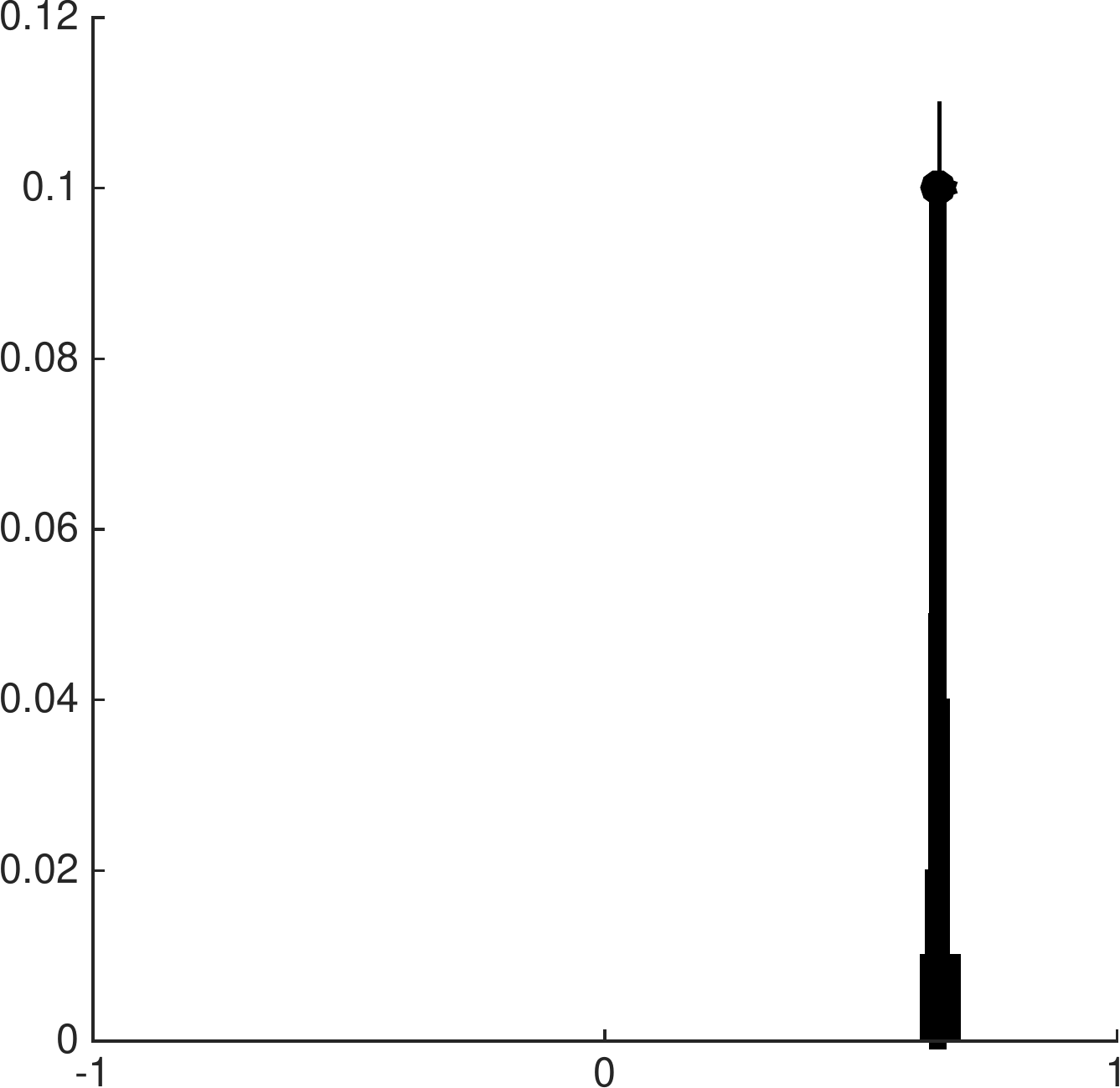}
}
\subfigure[Turb.~Intensity]{
\includegraphics[width=0.22\textwidth]{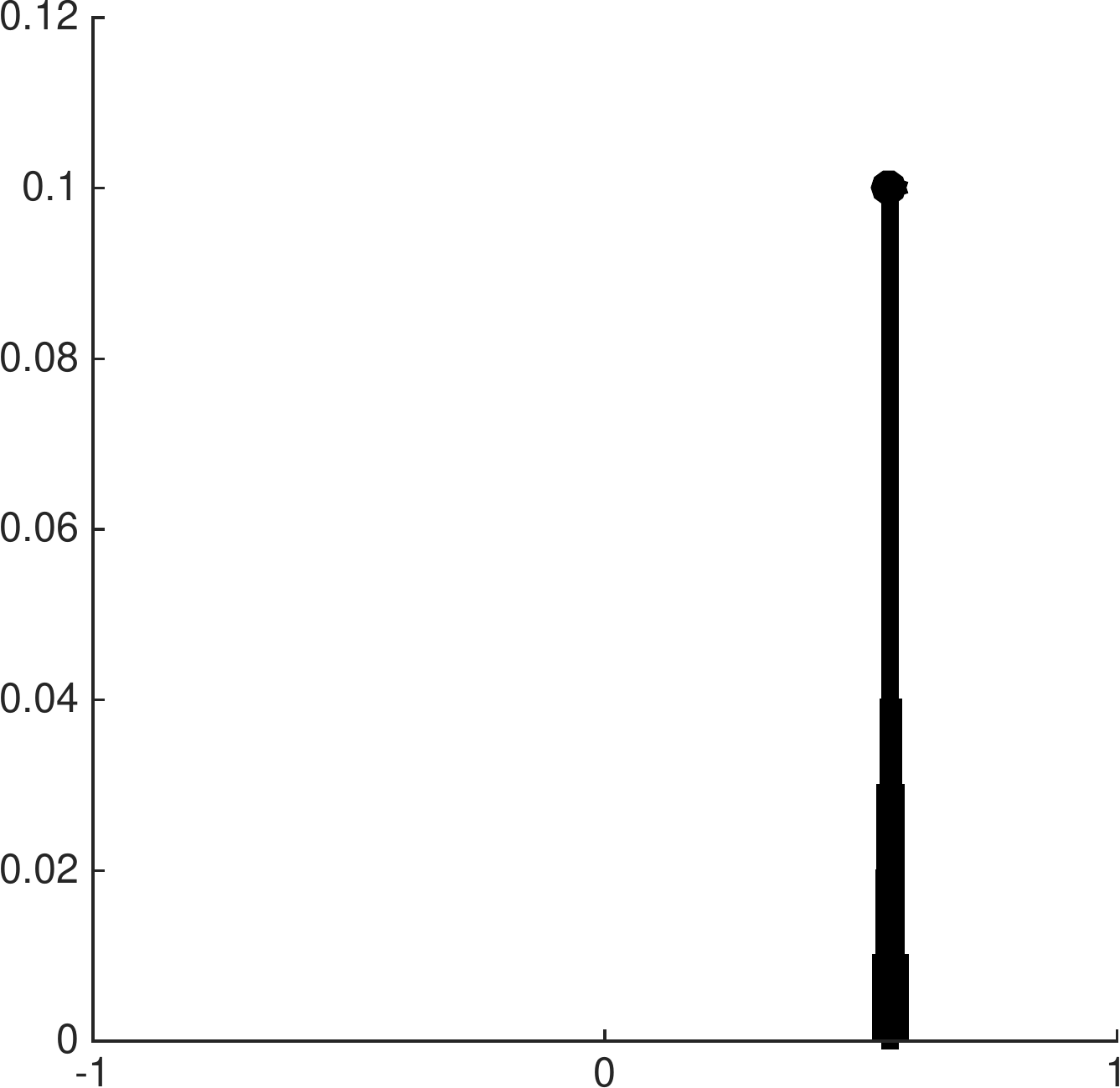}
}
\subfigure[Turb.~Len.~Scale]{
\includegraphics[width=0.22\textwidth]{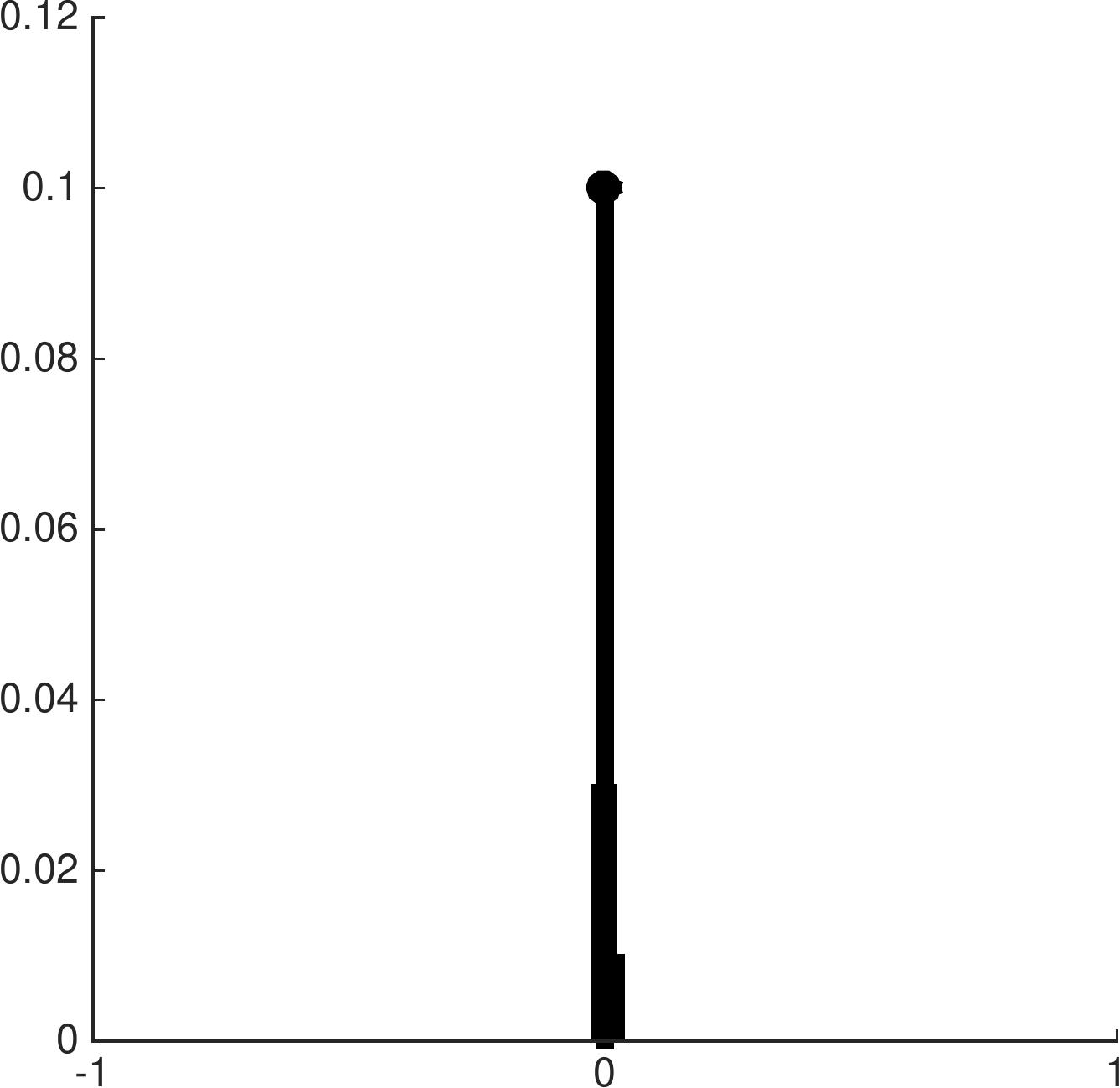}
}
\subfigure[Stag.~Pres.]{
\includegraphics[width=0.22\textwidth]{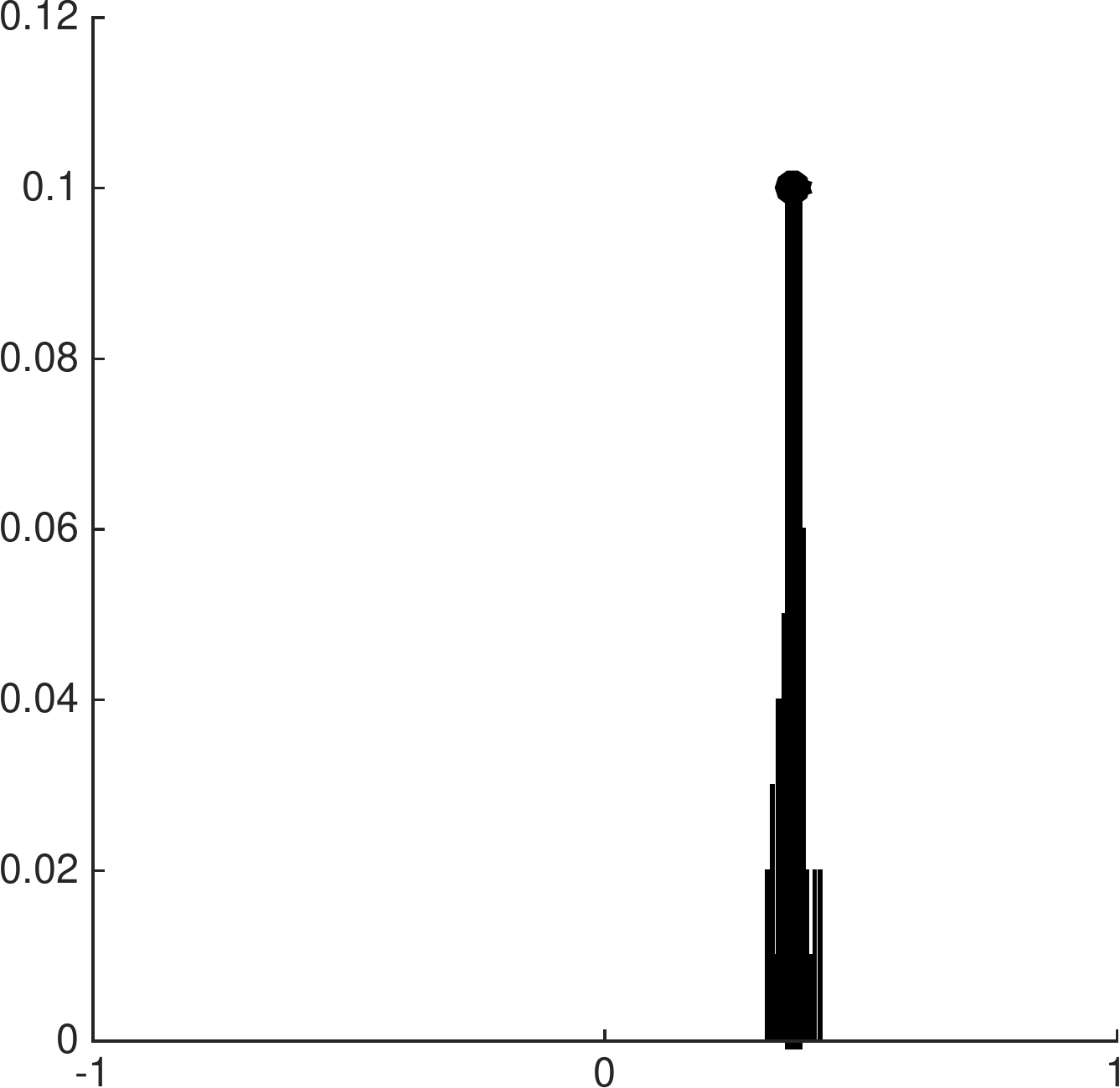}
}\\
\subfigure[Stag.~Enth.]{
\includegraphics[width=0.22\textwidth]{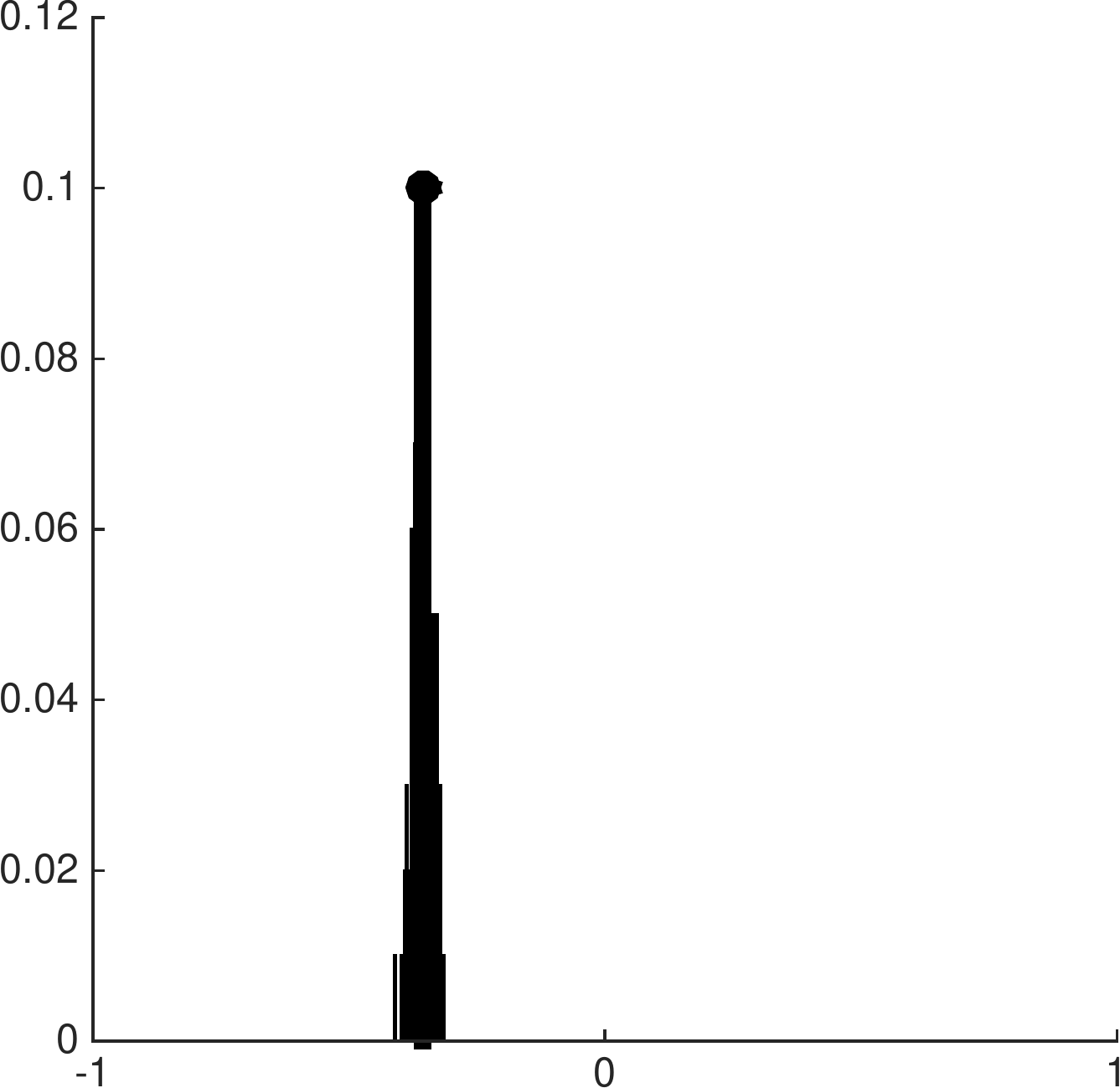}
}
\subfigure[Cowl Trans.]{
\includegraphics[width=0.22\textwidth]{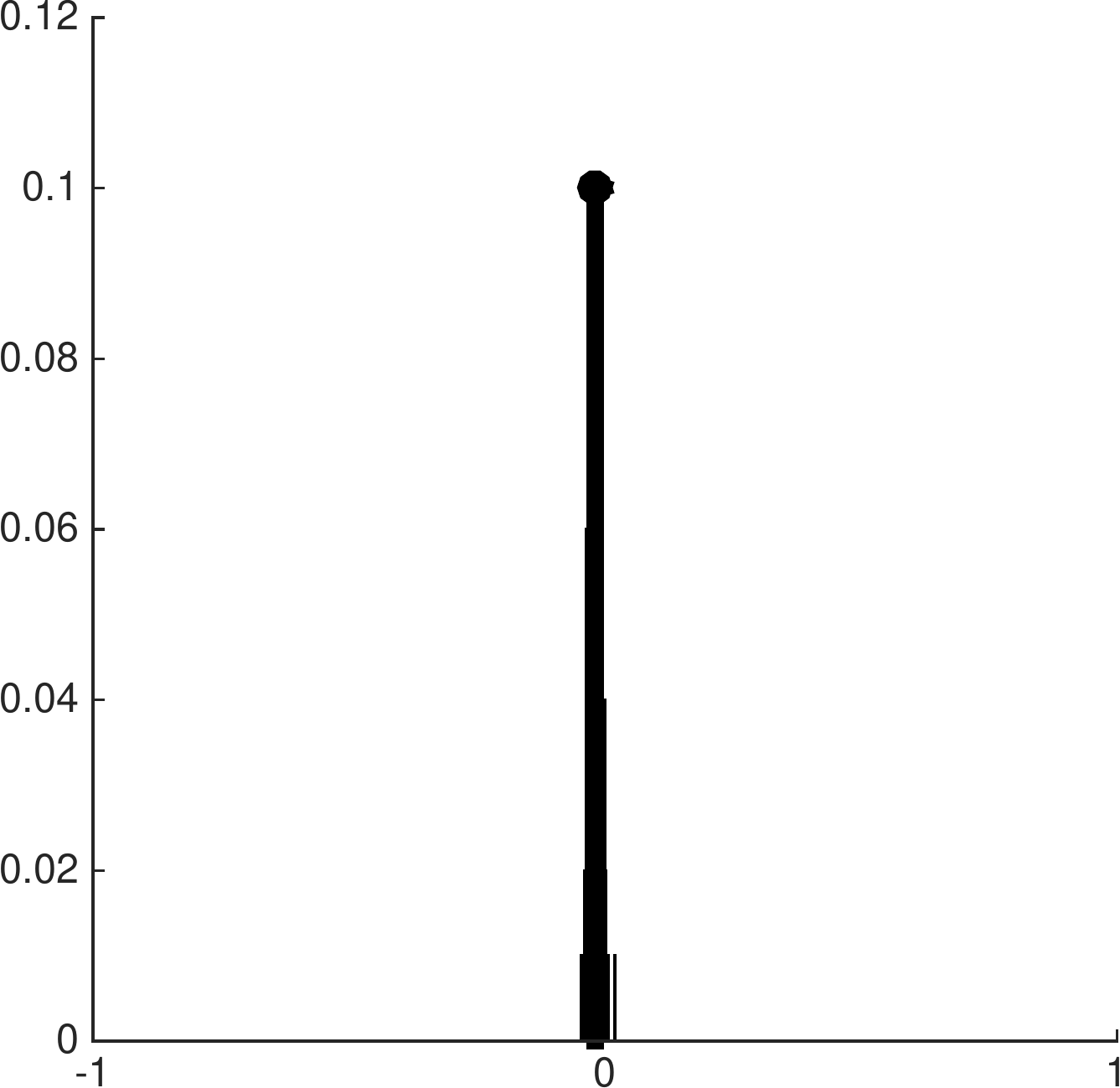}
}
\subfigure[Ramp Trans.]{
\includegraphics[width=0.22\textwidth]{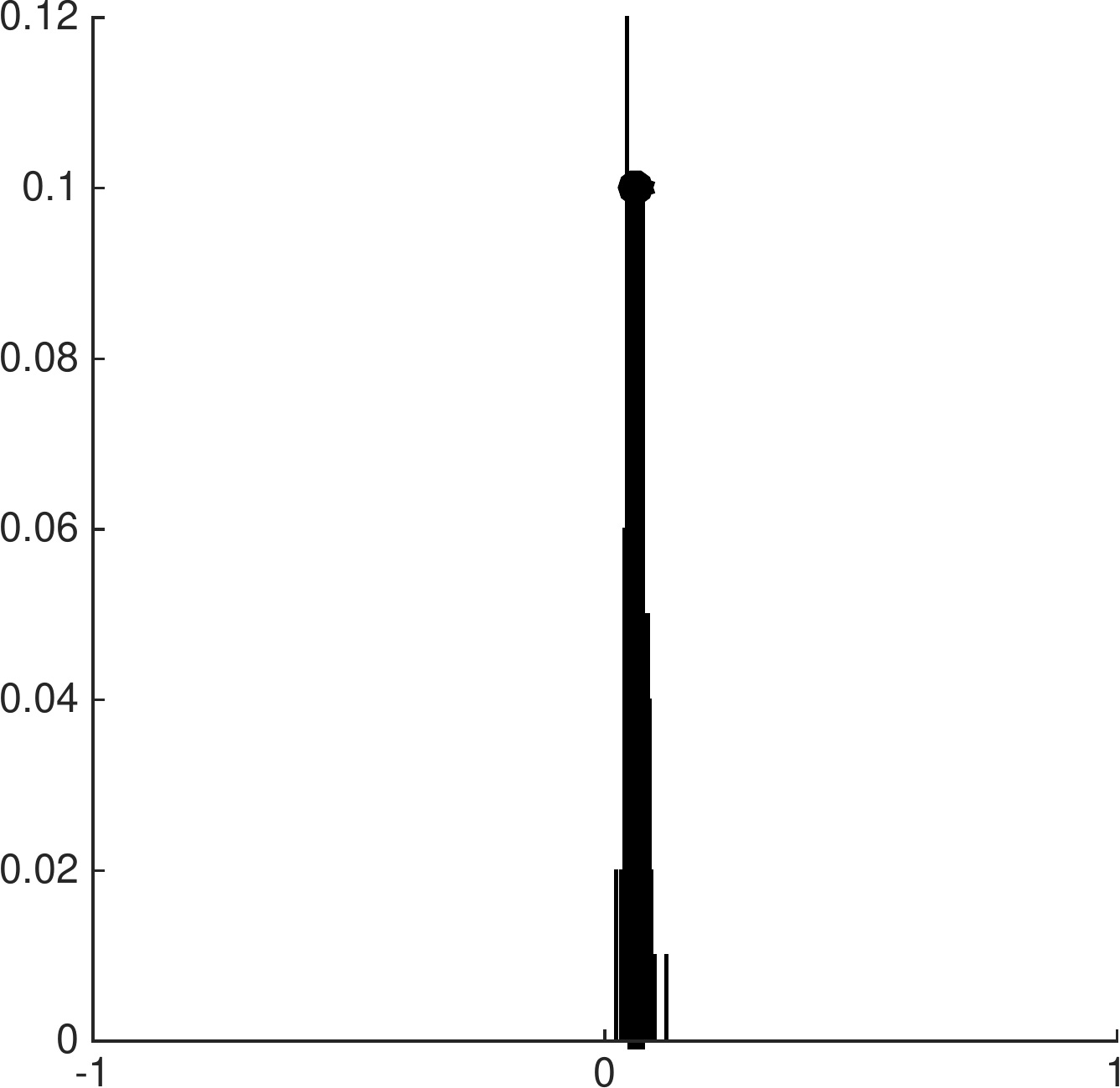}
}
\caption{Bootstrap histograms of the components of the active subspace vector $\vw$ for $P_{0,\rm H2}=4.8$ bar. The black stems are the computed components of $\vw$ from Table \ref{tab:dir}. The caption of each subfigure names the specific inflow parameter. The sharp peaks around each of the stems provides confidence that the computed $\vw$ are stable.}
\label{fig:bootstrap0}
\end{figure}

\begin{figure}[htb]
\centering
\subfigure[Angle of Attack]{
\includegraphics[width=0.22\textwidth]{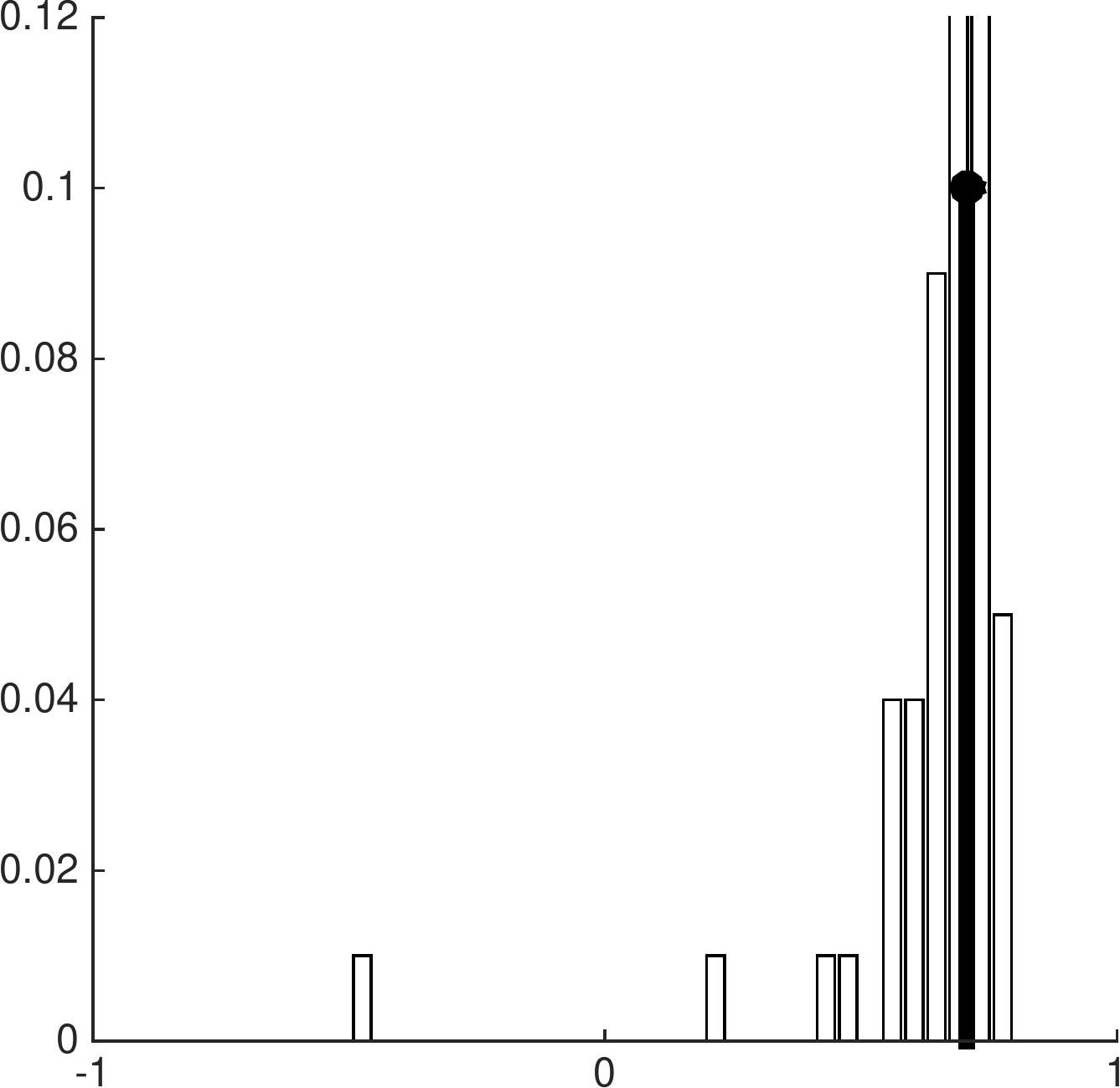}
}
\subfigure[Turb.~Intensity]{
\includegraphics[width=0.22\textwidth]{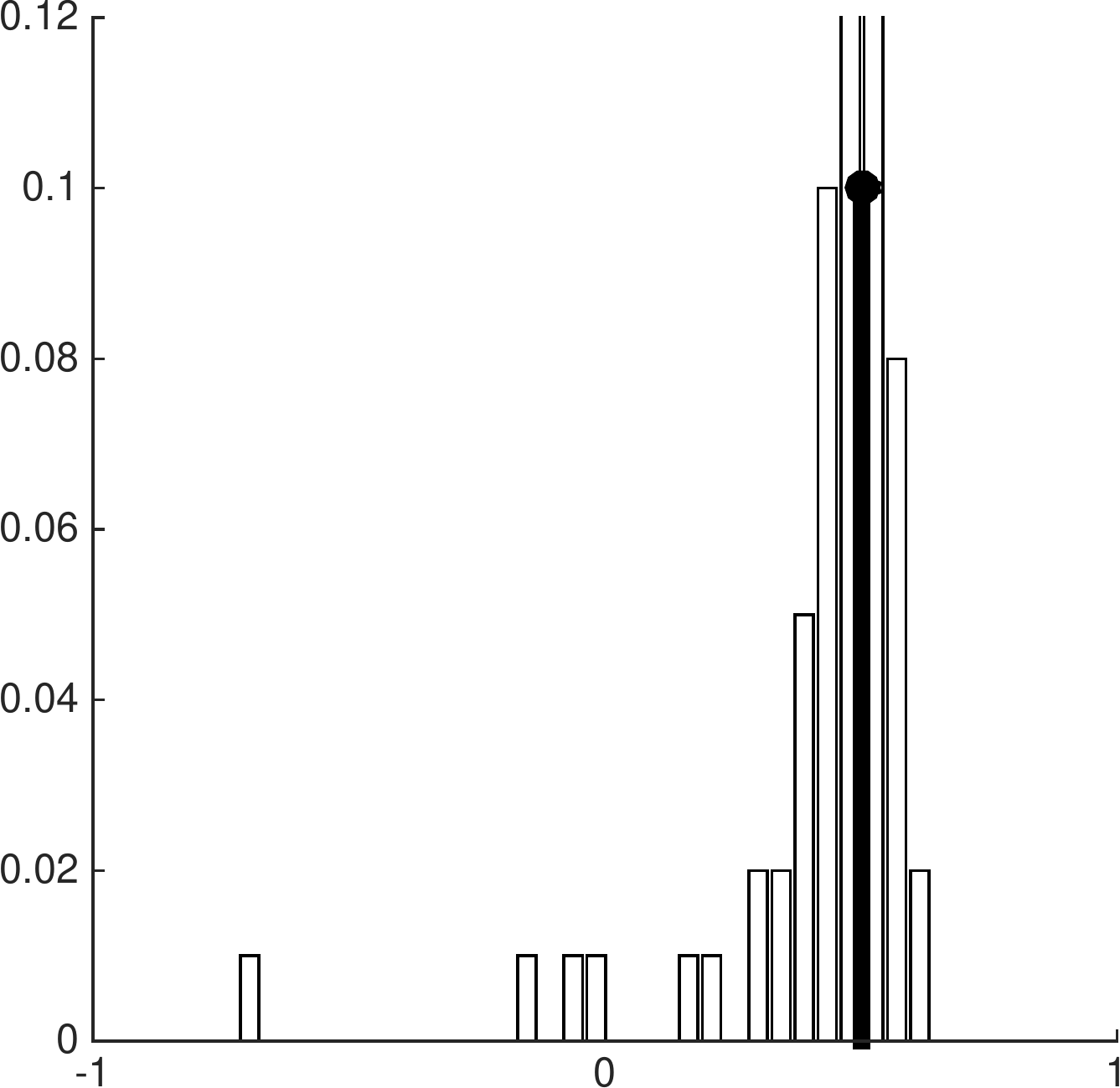}
}
\subfigure[Turb.~Len.~Scale]{
\includegraphics[width=0.22\textwidth]{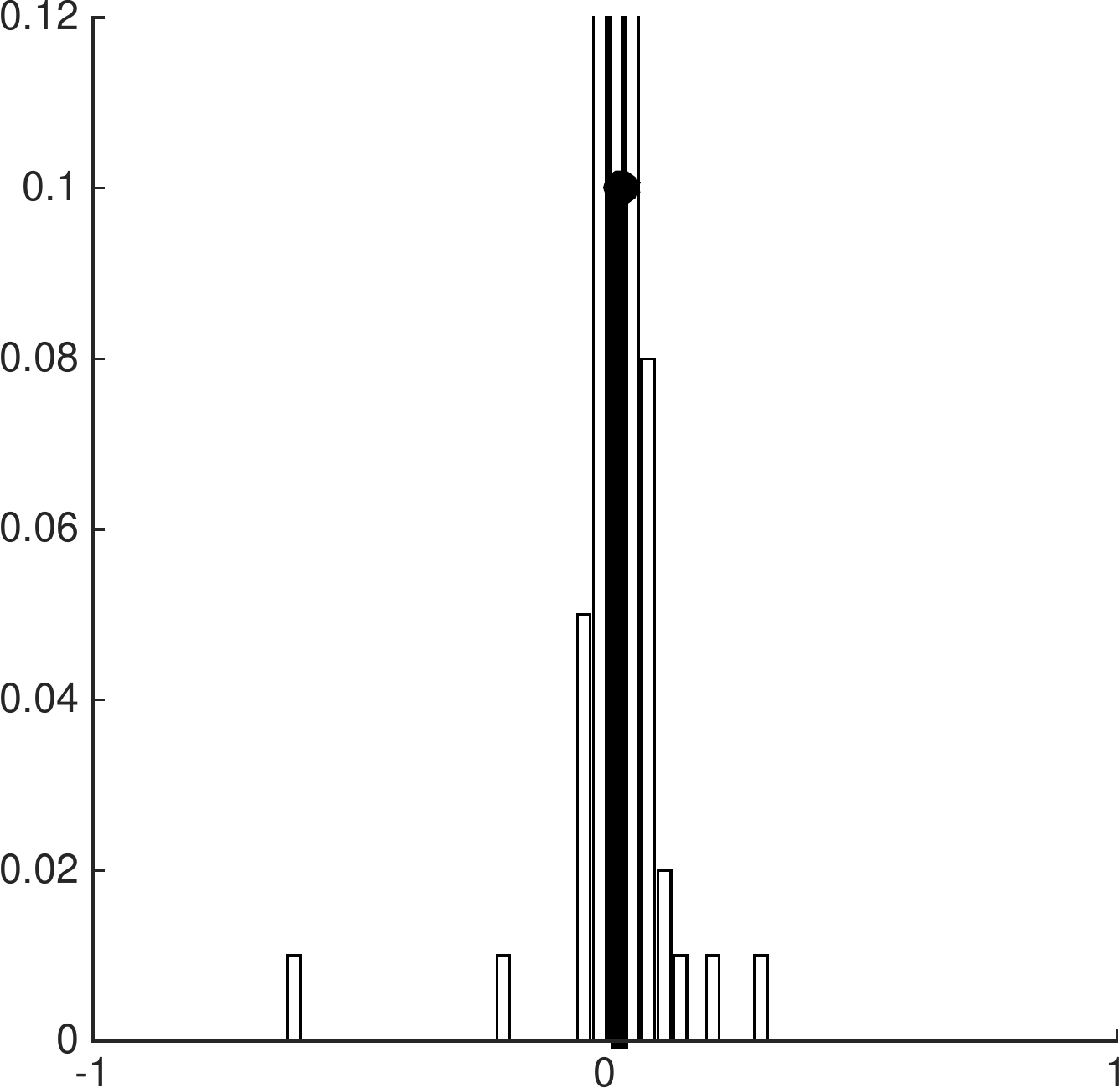}
}
\subfigure[Stag.~Pres.]{
\includegraphics[width=0.22\textwidth]{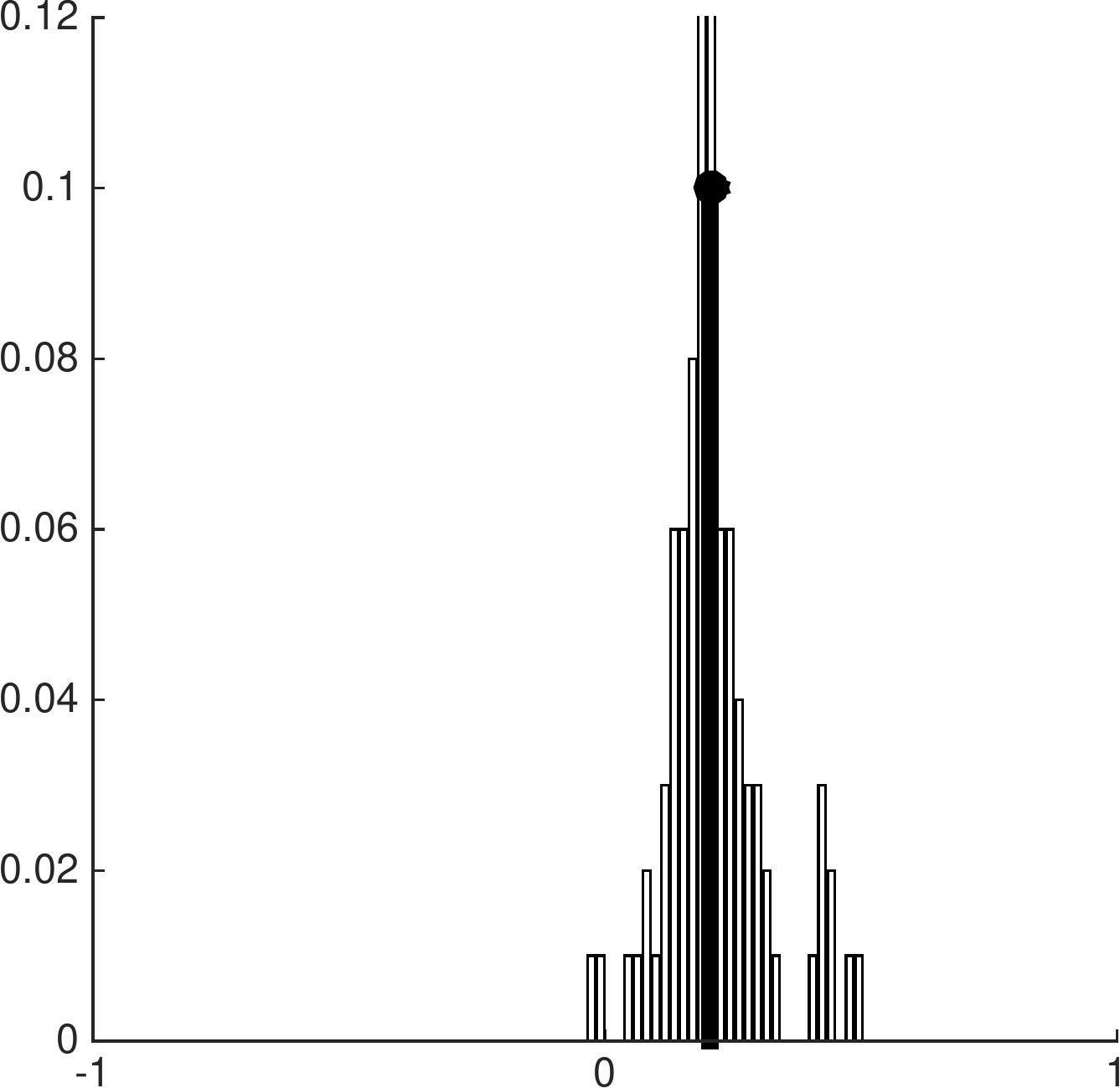}
}\\
\subfigure[Stag.~Enth.]{
\includegraphics[width=0.22\textwidth]{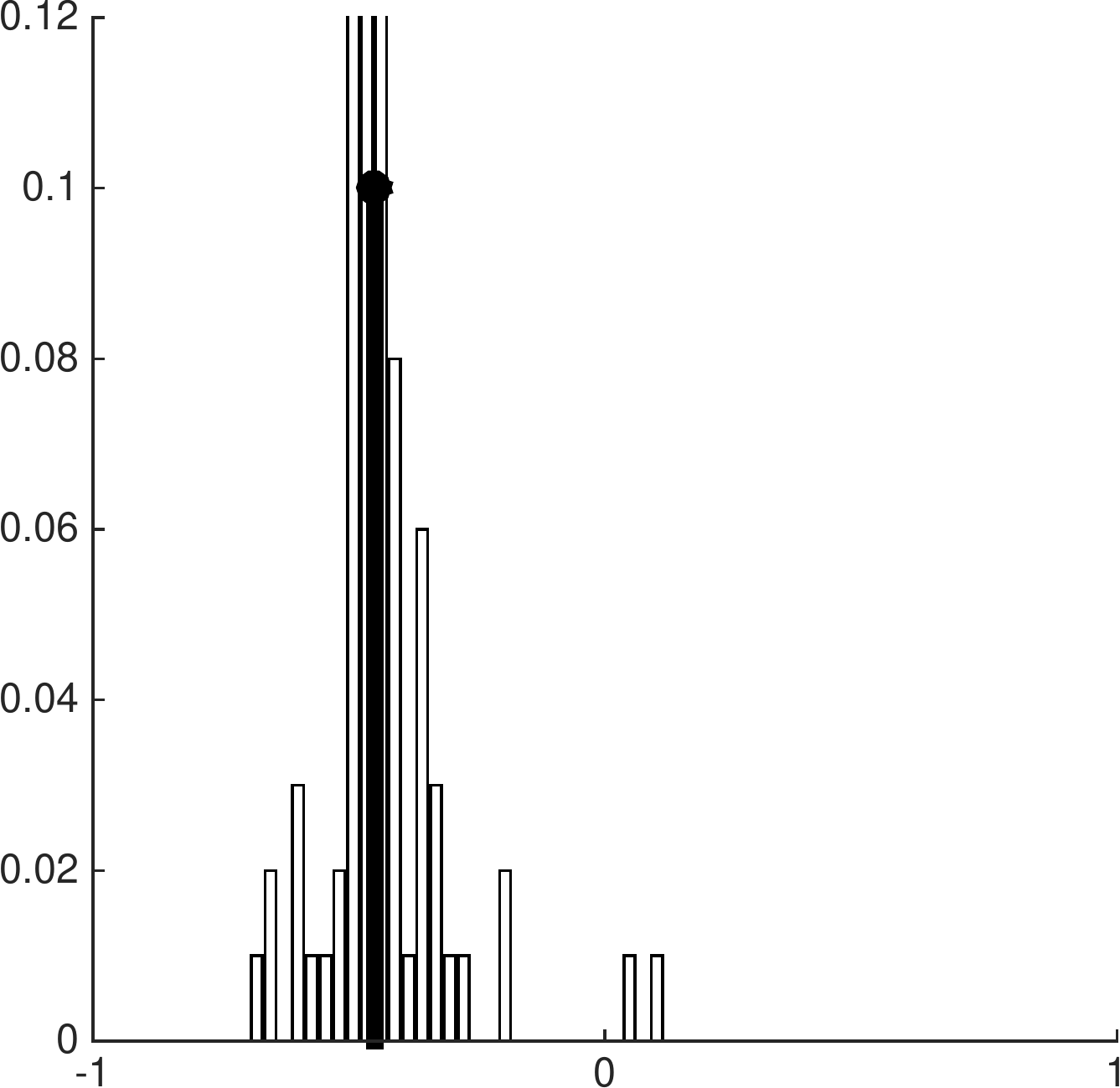}
}
\subfigure[Cowl Trans.]{
\includegraphics[width=0.22\textwidth]{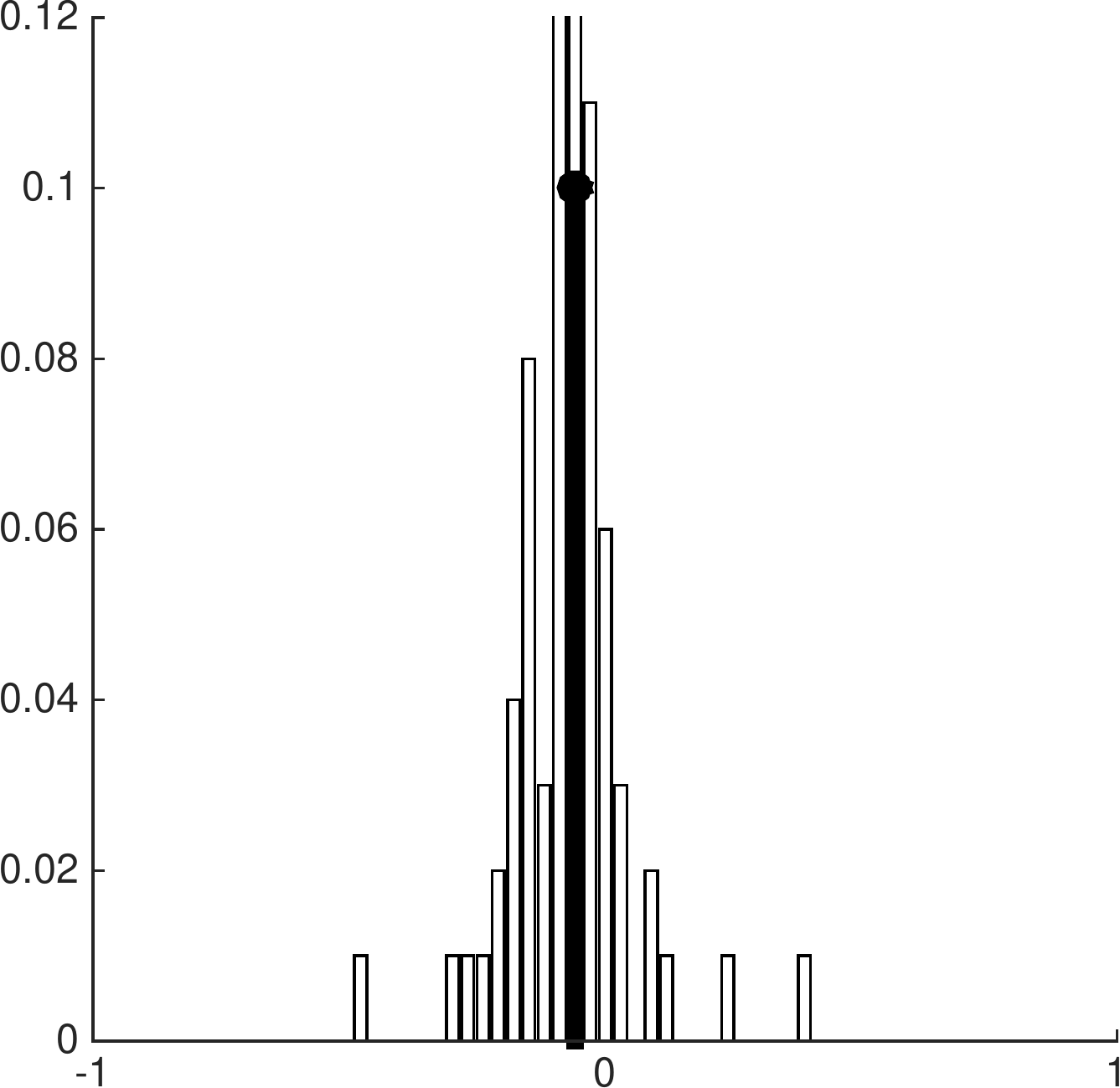}
}
\subfigure[Ramp Trans.]{
\includegraphics[width=0.22\textwidth]{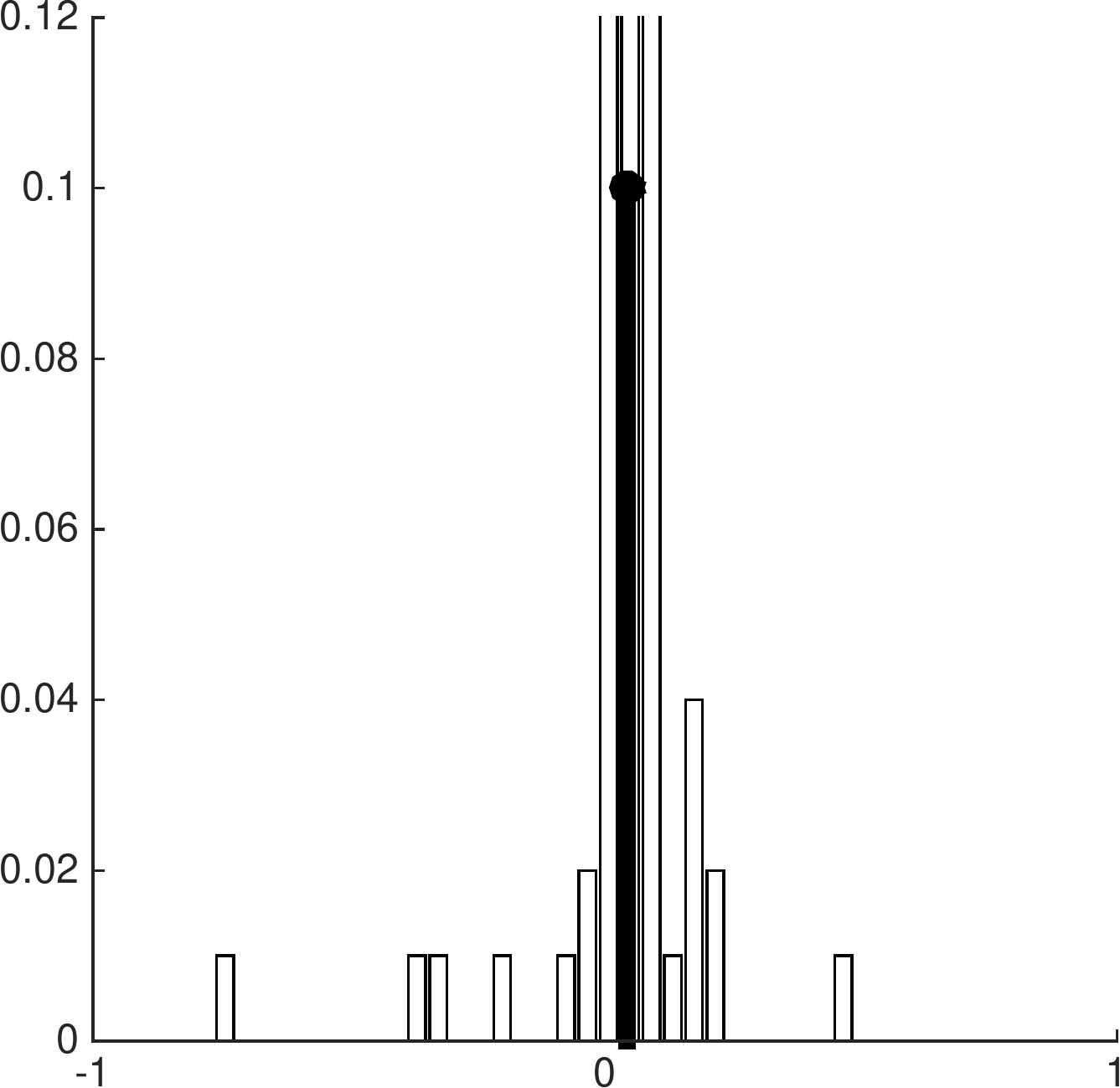}
}
\caption{Bootstrap histograms of the components of the active subspace vector $\vw$ for $P_{0,\rm H2}=5.6$ bar. The black stems are the computed components of $\vw$ from Table \ref{tab:dir}. The caption of each subfigure names the specific inflow parameter. The sharp peaks around each of the stems provides confidence that the computed $\vw$ are stable. The relatively large spread compared to Figure \ref{fig:bootstrap0} is due to the relatively small number ($M=14$) of samples.}
\label{fig:bootstrap1}
\end{figure}

\noindent We apply the bootstrap procedure from Algorithm \ref{alg:bs} to $\vw$ for the HyShot II exit pressure using $N=100$ bootstrap replicates. Note that the bootstrap uses only the available simulation runs, so its cost is negligible. The bootstrap histograms along with stem plots of the components of $\vw$ from Table \ref{tab:dir} are shown for $P_{0,\rm H2}=4.8$ bar in Figure \ref{fig:bootstrap0} and for $P_{0,\rm H2}=5.6$ bar in Figure \ref{fig:bootstrap1}. The histograms' sharp peaks around the stems suggests confidence in the computed directions. The relatively large ranges in Figure \ref{fig:bootstrap1} are a result of the low oversampling factor ($M=2m$) used to fit the linear approximation; with too few samples to draw from, the probability is greater that a bootstrap replicate is not representative of the true values of $\vw$. 

We can study how bootstrap variability in $\vw$ affects the perceived relationship in the summary plot. We use the 100 bootstrap replicates $\vw_k$ to plot the exit pressures $f_j$ against the corresponding active variables $\vw_k^T\vx_j$. The result is a horizontal scatter of gray dots around each point in Figure \ref{fig:as}. The scatter provides a visual indication of how $\vw$'s variability affects the perceived relationship between the active variable $\vw^T\vx$ and the exit pressure $f$. There is hardly any scatter in the case of $P_{0,\rm H2}=4.8$ bar with $M=50$ samples; the 100 gray circles per original sample (black circle) are barely visible behind the original samples. For $P_{0,\rm H2}=5.6$ bar, the large spread relative to the range of the active variable is due to the low sampling ($M=2m$) used to fit the linear approximation. 

In the lower fuel pressure case, the bootstrap instills confidence as intended, because 50 samples is sufficient to produce meaningful results with resampling. In the high fuel pressure case, 14 samples is not enough for the bootstrap to instill confidence---although 14 samples in the summary plot (Figure \ref{fig:as2}) still show a strong univariate relationship between the active variable and the exit pressure. Despite the bootstrap variability in the second case, three things provide confidence of the one-dimensional active subspace: (i) the strong validation in the first case with 50 runs and low bootstrap variability, (ii) the physical intuition that the relationship between inputs and outputs should not change dramatically between the two cases, and (iii) the apparent relationship in the summary plot Figure \ref{fig:as2}.

%


We exploit the one-dimensional active subspace revealed in the summary plots (Figure \ref{fig:as}) and confirmed by the bootstrap to quantify uncertainty in the HyShot II's simulated prediction of exit pressure. The uncertainty quantification is comprised of three computations: (i) estimating the range of possible exit pressures, (ii) identifying safe operating conditions, and (iii) estimating a cumulative distribution function for the exit pressure. 

\subsection{Approximating the range of exit pressures}
\label{sec:bounds}
\noindent The conditions noted at the beginning of Section \ref{sec:uq} are common in design optimization problems with expensive computer simulations. Several approaches attempt to optimize with the aid of a response surface~\cite{billups2013,wild2008orbit,jones2001taxonomy}. Unfortunately, our simulation is too expensive to use the iterative procedures that refine in possible regions of the optima. 

The summary plots in Figure \ref{fig:as} suggest that the exit pressure $f(\vx)$ is a monotonic function of the active variable $\vw^T\vx$. A univariate, monotonic, and continuous function defined on a bounded interval can be bounded by the function values at the interval boundaries. We use the following heuristic to estimate the maximum and minimum exit pressures. Define
\begin{equation}
\label{eq:xmax}
\vx_{\text{max}} \;=\; \argmax{-1\leq \vx \leq 1} \; \vw^T\vx,
\qquad
f_{\text{max}} \;=\; f(\vx_{\text{max}}).
\end{equation}
Similarly define
\begin{equation}
\label{eq:xmin}
\vx_{\text{min}} \;=\; \argmin{-1\leq \vx \leq 1} \; \vw^T\vx,
\qquad
f_{\text{min}} \;=\; f(\vx_{\text{min}}).
\end{equation}
The components of $\vx_{\text{max}}$ and $\vx_{\text{min}}$ are either -1 or 1, and the two points correspond to opposite corners of the hypercube, i.e., $\vx_{\text{min}}=-\vx_{\text{max}}$. The signs of $\vx_{\text{max}}$'s components are determined entirely by the corresponding signs of $\vw$'s components. 

The range $[f_{\text{min}},f_{\text{max}}]$ provides our estimate for the range of exit pressures from the HyShot II model. Computing this interval requires running the model four times beyond the initial set of samples used to compute $\vw$---two for each $P_{0,\rm H2}$---which is much cheaper than adaptively constructing a response surface. The values of $f_{\text{min}}$ and $f_{\text{max}}$ are shown in Figure \ref{fig:as} as black squares along with the samples (black circles) used to compute $\vw$ defining the active subspace---all plotted against the active variable. These runs confirm the monotonic structure perceived in the summary plot. At worst, these runs bound the initial exit pressure samples. At best, they provide estimates of the range of exit pressures over all values of the input parameters. Admittedly, this approach is not far from \emph{guess-and-check}. But the guesses were informed by the structure revealed in the summary plot. Checking the necessary conditions for stationarity of these points is not feasible since we do not have first or second derivatives. One can interpret this heuristic as one step of a least-squares-fit linear approximation-based approach for derivative-free optimization with one global linear model~\cite{conn2009introduction}. However, the important difference is that the summary plot provides evidence that the one global linear model is sufficient to estimate the optima. 

It is natural to wonder if the estimated range $[f_{\text{min}},f_{\text{max}}]$ is strongly affected by the randomness in $\vw$ originating in the random sampling in Algorithm \ref{alg:as}. If some of $\vw$'s components are sufficiently far from zero, then small perturbations in $\vw$ do not change the signs of the corresponding components of  $\vx_{\text{min}}$ and $\vx_{\text{max}}$. If a component of $\vw$ is small in magnitude---such as components 2, 6, and 7 in Table \ref{tab:dir}---then small, random changes may change the corresponding components of $\vx_{\text{min}}$ and $\vx_{\text{max}}$ from -1 to 1 or vice versa. However, by the active subspace's construction, changes in these parameters do not change the exit pressure as much as changes in parameters with large $\vw$ components. This gives confidence that the estimated range is stable under small changes in $\vw$. 

\subsection{Constraining the exit pressure}
\noindent Consider the following exercise in safety engineering. Suppose that the scramjet operates safely when the exit pressure is below 2.8 bar, but it nears unsafe operation above 2.8 bar. With the one-dimensional active subspace, we can quickly characterize the parameter regime that produces exit pressures below the threshold of 2.8 bars. 

The first step is to build a response surface model of the exit pressure as a function of the active variable $\vw^T\vx$. We could try to construct a response surface of all seven input parameters. But with only 52 model runs for $P_{0,\rm H2}=4.8$ bar (50 runs to compute $\vw$ and 2 runs to estimate the range of $f$) and 16 model runs for $P_{0,\rm H2}=5.6$ bar (14 for $\vw$ and 2 for the estimated range), our modeling choices would be very limited. The apparent low departure from a univariate relationship in Figure \ref{fig:as} suggests that we can construct a useful response of just the active variable as in \eqref{eq:1d}. In particular, we can use the exit pressure samples $f_j$ and corresponding active variable values $\vw^T\vx_j$ to model the exit pressure as a univariate quadratic polynomial of the active variable. 

In statistical regression, the coefficient of determination (i.e., $R^2$) is one of several metrics for the quality of the model~\cite{weisberg2005applied}. This coefficient is interpreted as the proportion of variance in the response (the output) that comes from variance in the predictors (the inputs); $R^2$ values near 1 indicate a good model for the data. This interpretation is not appropriate in our case, since there is no random noise in the $f_j$'s. Nevertheless, we report the $R^2$ values for the quadratic model of the active variable and treat them as a measure of discrepancy between the model and the data: $R^2=0.993$ for $P_{0,\rm H2}=4.8$ bar and $R^2=0.998$ for $P_{0,\rm H2}=5.6$ bar. 

\begin{figure}[htb]
\centering
\subfigure[$P_{0,\rm H2}=4.8$ bar]{
\includegraphics[width=0.46\textwidth]{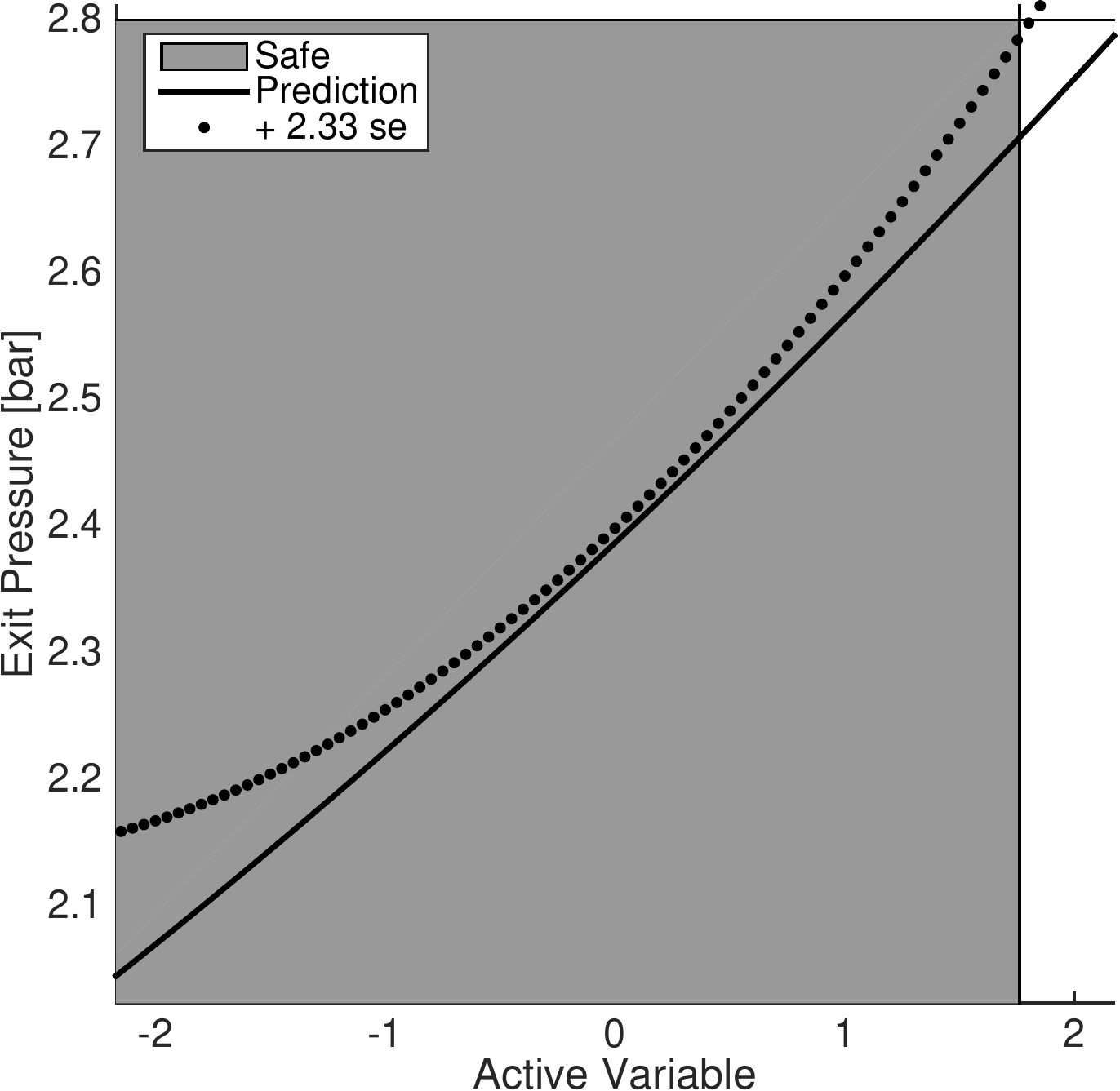}
}
\subfigure[$P_{0,\rm H2}=5.6$ bar]{
\includegraphics[width=0.46\textwidth]{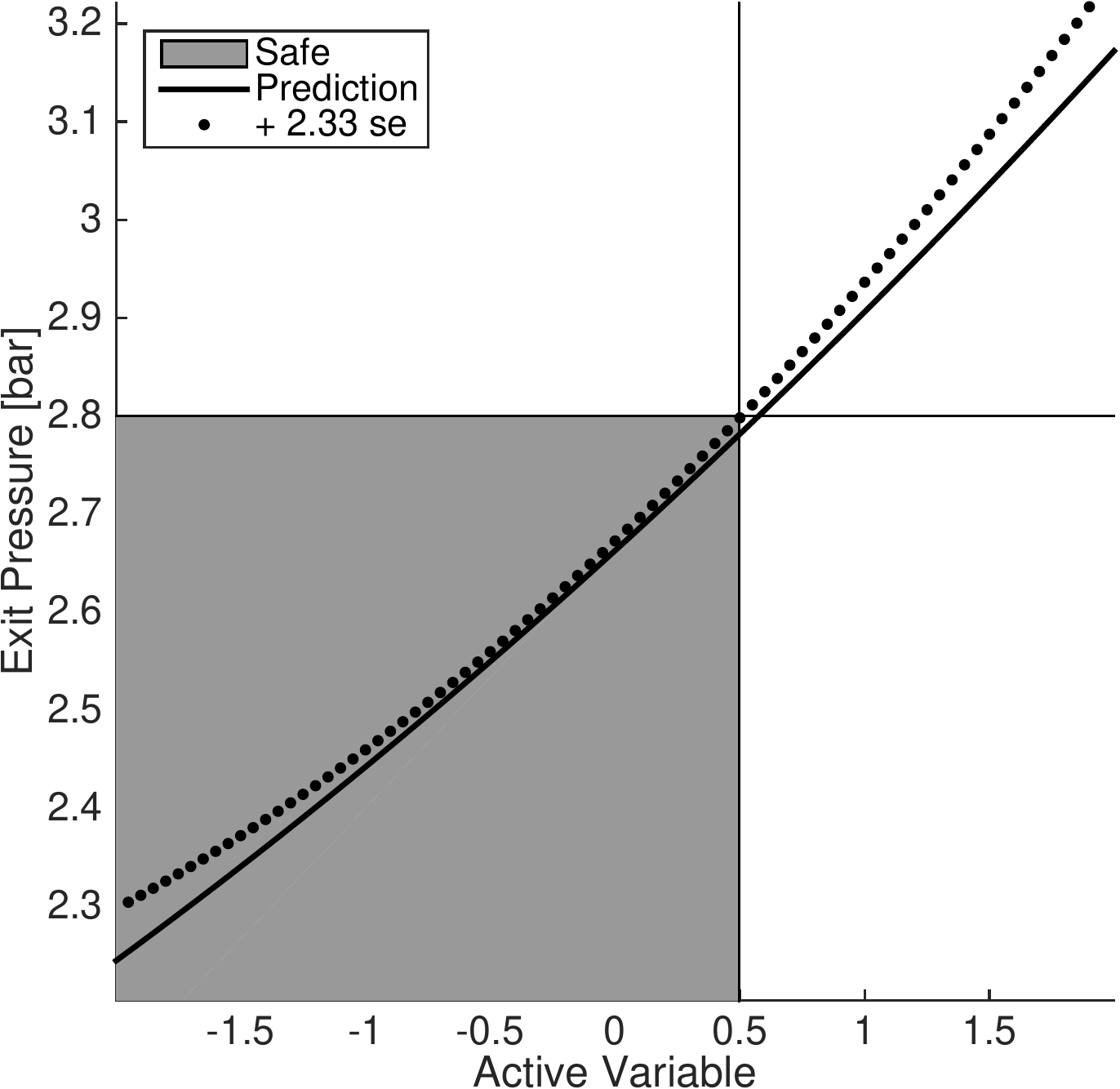}
}
\caption{A quadratic polynomial models the relationship between the active variable $\vw^T\vx$ and the exit pressure. The solid line shows the quadratic model. The dotted line shows the computed upper 99\% confidence bound. We find the value of the active variable $\ymax$ where the upper confidence limit crosses the safety threshold of 2.8 bars. All values of the active variable less than $\ymax$ produce pressures within the safety limit. The set of safe input variables is shown in \eqref{eq:safeset}.}
\label{fig:safe}
\end{figure}

Regression predictions often include frequentist confidence bounds. We can compute the upper 99\% confidence bound and treat it as a conservative factor when seeking the safe operating conditions.\footnote{Since there is no noise in the $f_j$'s, the statistical interpretation of the upper confidence bound is not valid; for example, we cannot interpret the 99\% confidence bounds as a random interval that contains the exit pressure with probability 0.99.} We find the largest value of the active variable such that the upper confidence bound from the quadratic approximation is less than the threshold of 2.8 bars. Figure \ref{fig:safe} shows the quadratic approximation and its upper 99\% confidence limit for both values of $P_{0,\rm H2}$. The shaded region identifies the values of the active variable that produce exit pressures at or below the pressure threshold. 

Let $\ymax$ be the value of the active variable where the upper confidence bound crosses the pressure threshold. Then the safe region of the normalized input parameters is the set $\sS$ defined as 
\begin{equation}
\label{eq:safeset}
\sS \;=\; \{\vx \,:\, \vw^T\vx\leq \ymax,\; -1\leq \vx \leq 1\}.
\end{equation}
One can easily shift and scale this region to the space of the HyShot II model's input space for a physical interpretation. The linear inequality constraint implies that the inputs are related with respect to the exit pressure. For example, the range of safe angles of attack depends on the other input variables. The presence of the active subspace and the quality of the quadratic approximation enable us to simply characterize these relationships. 

The safe set $\sS$ defined in \eqref{eq:safeset} is like a seven-dimensional box with the top chopped off by the hyperplane $\vx^T\vw\leq\ymax$ (though the notion of ``top'' is problematic in seven dimensions). We can identify a new set of independent ranges for the input variables such that all inputs within those ranges produce exit pressures below the 2.8 bar safety threshold---according to the quadratic response surface model. This is like finding the largest seven-dimensional box that fits inside the set $\sS$. More precisely, we solve the following optimization problem,
\begin{equation}
\begin{array}{lc}
\underset{\vx}{\mathrm{maximize}} & \prod_{i=1}^m |x_i-x_{i,\text{min}}|,\\
\st & \vx\in\sS,
\end{array}
\end{equation}
where $x_{i,\text{min}}$ are the components of the minimizer $\vx_{\text{min}}$ from \eqref{eq:xmin}. The maximizing components define the corner of the largest hyperrectangle opposite the corner $\vx_{\text{min}}$. We can shift and scale this hyperrectangle to the space of model parameters for physical interpretation. Table \ref{tab:safe_ranges} shows the ranges of parameters affected by the safety condition. 

\begin{table}[h!tb]
\begin{center}
\begin{tabular}{l | l | c c | c}
 & Parameter & Min & Max & Units \\ 
\hline
\multirow{2}{*}{$P_{0,\rm H2}=4.8$ bar} & Angle of Attack & 2.6 & 4.29 & deg. \\
 & Turbulence Intensity & 0.001 & 0.0188 & $\cdot$ \\
\hline
\multirow{3}{*}{$P_{0,\rm H2}=5.6$ bar} & Stagnation Enthalpy & 3.15 & 3.4280 & $\mbox{MJ}/\mbox{kg}$ \\
 & Angle of Attack & 2.6 & 3.54 & deg. \\
 & Turbulence Intensity & 0.001 & 0.013 & $\cdot$ \\
\end{tabular}
\end{center}
\caption{Updated ranges for input parameters that lead to safe scramjet operation with exit pressure, modeled with a quadratic function of the active variable, less than 2.8 bars. The parameters listed for each $P_{0,\rm H2}$ case are the parameters affected by the restriction in \eqref{eq:safeset}. Compare these ranges to those in Table \ref{tab:inflow_uq_ranges}.}
\label{tab:safe_ranges}
\end{table}

We can interpret such analysis as backward uncertainty propagation that characterizes safe inputs given a characterization of a safe output under the constraint that the inputs be independent. These spaces differ between the two fuel plenum pressures; see Figure \ref{fig:safe} and Table \ref{tab:safe_ranges}. For $P_{0,\rm H2}=4.8$ bar, only angle of attack and turbulence intensity are affected by the safety constraint on the pressure; angle of attack must be less than 4.29 degrees, and turbulence intensity must be less than 0.0188. For $P_{0,\rm H2}=5.6$ bar, the same parameters are constrained---angle of attack less than 3.54 and turbulent intensity less than 0.013---and a stricter minimum on stagnation enthalpy (greater than 3.15 MJ/Kg) appears. To explain this, we note that low free-stream stagnation temperature leads to conditions that are further away from the adiabatic frame temperature, which leads to increased heat deposition in the air-stream. This in turns leads to more deceleration and larger pressure increases in the chamber \cite{wendt1997}.

\subsection{Cumulative distribution function}
\label{sec:cdf}

\begin{figure}[htb]
\centering
\subfigure[$P_{0,\rm H2}=4.8$ bar]{
\includegraphics[width=0.46\textwidth]{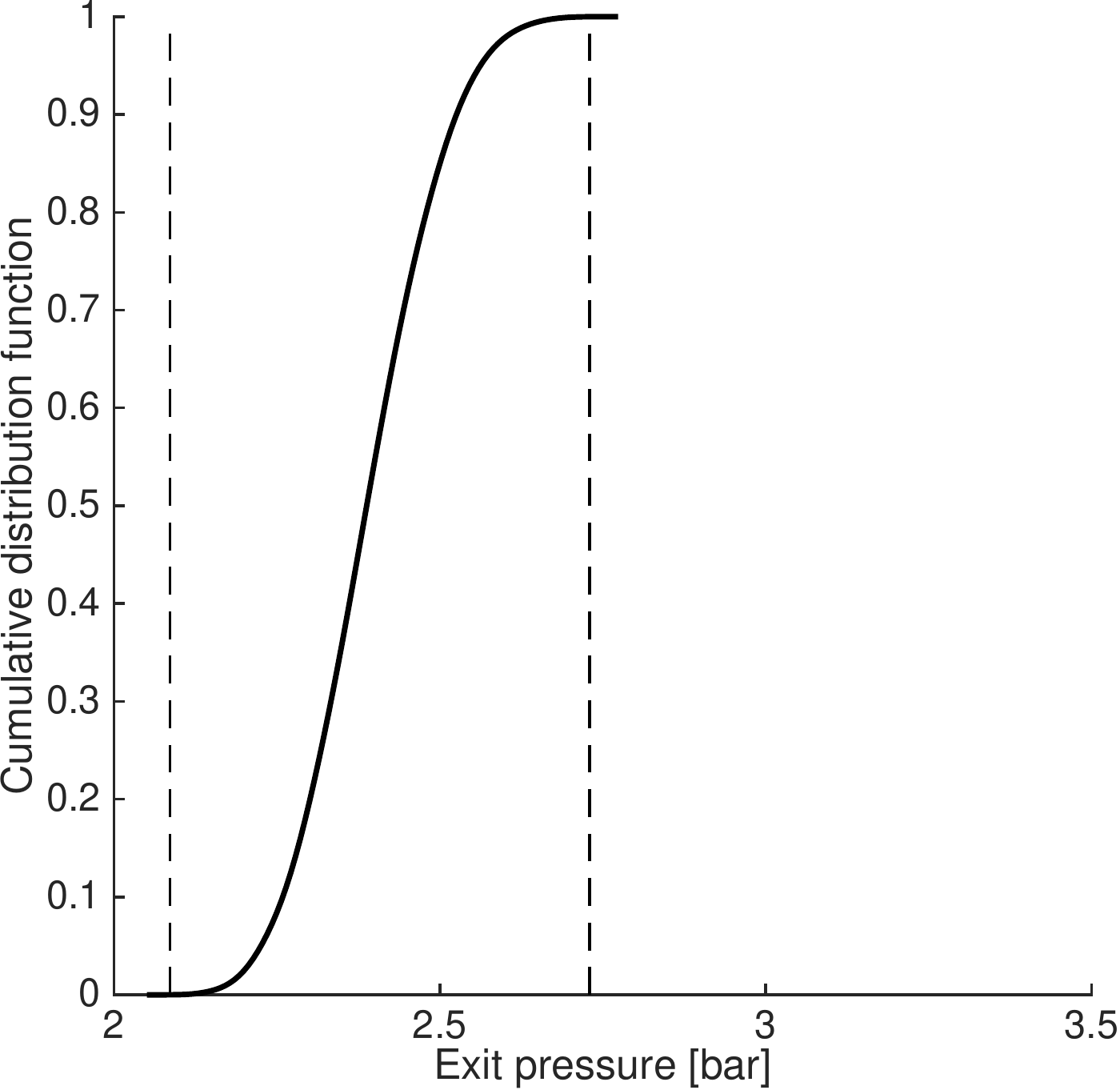}
}
\subfigure[$P_{0,\rm H2}=5.6$ bar]{
\includegraphics[width=0.46\textwidth]{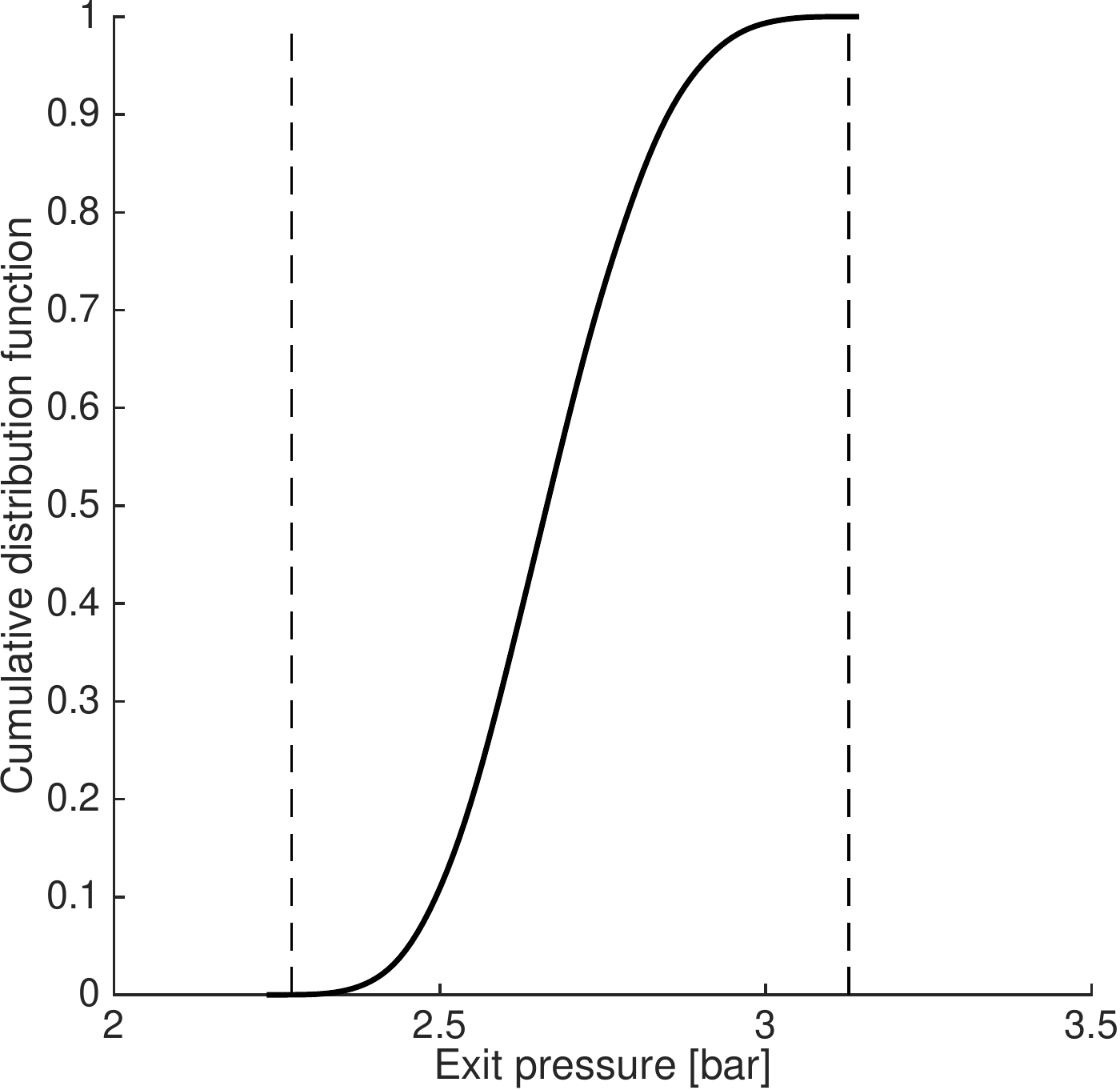}
}
\caption{Estimated cumulative distribution functions for exit pressure at both values of fuel plenum pressure $P_{0,\rm H2}$. These are estimated with Gaussian kernel density estimates; the samples are drawn from the quadratic approximation of exit pressure as a function of the active variable. Vertical lines show the estimated bounds for each case estimated in Section \ref{sec:bounds}.}
\label{fig:cdf}
\end{figure}

\noindent We use the quadratic approximation of the active variable to estimate a cumulative distribution function of the exit pressure for both fuel plenum pressures. We draw 5000 samples independently and uniformly from the scramjet's seven-dimensional input space (i.e., according the the uniform density $\rho(\vx)$), and for each sample we evaluate (i) the active variable and (ii) the quadratic approximation of exit pressure as a function of the active variable. Figure \ref{fig:cdf} shows estimates of the cumulative distribution function computed from these samples with a Gaussian kernel density estimator. The vertical lines show the estimated upper and lower bounds computed in Section \ref{sec:bounds}. 

%

\section{Summary and discussion}
\label{sec:conclusions}

\noindent We present a numerical investigation of the reactive flow within a hydrogen-fueled scramjet with the objective of studying the effect of uncertainties in operating conditions on the overall performance. We carry out three-dimensional RANS simulations with a flamelet-based combustion model at two different fuel plenum pressures and record the pressure at the engine exit as a measure of system performance. We consider seven uncertain parameters; their ranges are justified by experimental evidence and/or expert opinions. 

The uncertainty quantification starts by identifying a one-dimensional reparameterization of the map between simulation parameters and the exit pressure. This reparameterization is based on a one-dimensional active subspace; the coefficients of a global, least-squares-fit linear approximation define the active subspace and are validated by both a summary plot and a bootstrap. Fitting an accurate linear approximation for a function of $m$ variables takes $\mathcal{O}(m)$ simulations. For the smaller fuel plenum pressure, we use 50 simulations to estimate the vector defining the active subspace; for the larger fuel plenum pressure, we use only 14 simulations. In both cases, the summary plot shows strong evidence of near one-to-one map between the linear combination of the inputs (i.e., the active variable) and the exit pressure. However, with only 14 runs, the bootstrap results show much greater variability due to insufficient sampling. We treat the components of the vector defining the active subspace as sensitivity metrics to gain insight into the physics. 

With a one-dimensional approximation and the perceived monotonic relationship in the summary plot, the computations necessary to quantify uncertainty become much easier. The minimum and maximum values for a univariate, continuous, monotonic function reside at the endpoints of the interval domain; we exploit this fact to estimate the range of possible exit pressures over all parameter values for both fuel pressure cases. We can validate a univariate response surface with the summary plot to complement standard quality metrics (e.g., residual norms). We use the validated, univariate, monotonic response surface to identify safe operating conditions in the active variable, and these constraints translate to the original model parameters. Lastly, we use the validated response surface as a surrogate to estimate a cumulative distribution function of the exit pressure given variability in the input parameters. To quantify uncertainty for both fuel pressure cases, we used 68 full simulations: 64 to identify the respective one-dimensional active subspaces and 4 to estimate the respective exit pressure ranges.

Quantifying uncertainty in such a complex and expensive simulation with seven independent inputs using so few runs is remarkable. It was only possible because of the low-dimensional structure revealed by the active subspace. If the summary plots had not shown a univariate, monotonic relationship, then our simple heuristics for estimating the exit pressure range and validating a response surface would not be justified. Thus, our approaches are only appropriate in the presence of this type of structure. The innovation in active subspaces is the potential to cheaply identify such exploitable structure to make possible otherwise infeasible computations.


\section*{Acknowledgments}
\noindent This research was funded by the U.S. Department of Energy [National Nuclear Security Administration] under Award No. NA28614. Additional computational resources were made possible by the following award, MRI-R2: Acquisition of a Hybrid CPU/GPU and Visualization Cluster for Multidisciplinary Studies in Transport Physics with Uncertainty Quantification. This award is funded under the American Recovery and Reinvestment Act of 2009 (Public Law 111-5). The first author's work is supported by the U.S. Department of Energy Office of Science, Office of Advanced Scientific Computing Research, Applied Mathematics program under Award Number DE-SC-0011077.

\section*{References}
\bibliographystyle{elsarticle-harv}
\bibliography{inflowUQ}

\end{document}